\documentclass[11pt,twoside]{amsart}
\usepackage{amsmath, amsthm, amscd, amsfonts, amssymb, graphicx}
\usepackage[bookmarksnumbered, plainpages]{hyperref}
\usepackage{mathrsfs}
\UseRawInputEncoding
\newtheorem{thm}{Theorem}[section]

\newtheorem{lem}[thm]{Lemma}

\newtheorem*{Proof}{Proof}

\numberwithin{equation}{section}

\begin{document}

\title{Invariant Ricci collineations associated to the Bott connections on three-dimensional Lorentzian Lie groups}
\author{Yanli Wang}

\thanks{{\scriptsize
\hskip -0.4 true cm \textit{2010 Mathematics Subject Classification:}
53C40; 53C42.
 }}

\maketitle

\begin{abstract}
 In this paper, we determine all left-invariant Ricci collineations associated to the Bott connection with three distributions on three-dimensional Lorentzian Lie groups .
\end{abstract}

\vskip 0.2 true cm


\pagestyle{myheadings}
\markboth{\rightline {\scriptsize Wang}}
         {\leftline{\scriptsize Ricci collineations}}

\bigskip
\bigskip

\section{ Introduction}
\indent "Symmetry" can be regarded as a one-parameter group of diffeomorphism of spacetime, and it preserves mathematical and physical quantity. In 1969, based on the different geometric objects' invariance properties, Katzin defined them as those vector fields X, such that leave the various relevant geometric quantities and classified Ricci collineations, curvature collineations in \cite{Ka}. Collineations are symmetry properties of space-times, so Ricci collineations i.e. Ricci symmetry, which is defined by $(L_{\xi}Ric)=0$. Because Ricci Collineations' close connection with the energy-momentum tensor, inspired more mathematicians' research interesting. After that Carot in \cite{JC} described the general proper Ricci collineations in type B warped space-times for spherically symmetric metrics in 1997. Next year Qadir and Ziad discussed the Ricci collineations of spherically symmetric spacetimes in \cite{QZ}. Ricci Collineations have been discussed and determined in more different spacetimes and for various other ricci tensor in \cite{CN, AC, HT, HT1}.\\
\indent In \cite{F,J,RK}, the Bott connection has been introduced. And three-dimensional Lorentzian Lie groups had been divided into $\{G_i\}_{i=1,\cdots,7}$ in \cite{Ca1,CP}, by this classification, Yong Wang defined a product structure on three-dimensional Lorentzian Lie groups and compute canonical connections and Kobayashi-Nomizu connections and their curvature with this product structure in \cite{Wa}. He also classified algebraic Ricci solitons associated to canonical connections and Kobayashi-Nomizu connections on three-dimensional Lorentzian Lie groups with this product structure. Simially, in \cite{Wu} Tong Wu computed the Bott connections and their Ricci tensor and classified affine solitons associated to the Bott connections. In this paper, we determine all left-invariant Ricci collineations associated to three different Bott connections which is defined by three different distributions in \cite{Wu} on $\{G_i\}_{i=1,\cdots,7}$ ( $\{G_i\}_{i=1,\cdots,4}$ is three-dimensional unimodular Lorentzian Lie groups.$\{G_i\}_{i=5,\cdots,7}$ is three-dimensional non-unimodular Lorentzian Lie groups).\\
\indent In section 2, we recall the definition of the Bott connection $\nabla^{B_1}$ with the first distribution, and then give left-invariant Ricci collineations associated to the Bott connection $\nabla^{B_1}$ on three-dimensional Lorentzian unimodular and non-unimodular Lie groups. In section 3, we recall the definition of the Bott connection $\nabla^{B_2}$ with the second distribution, and then give left-invariant Ricci collineations associated to the Bott connection $\nabla^{B_2}$ on three-dimensional Lorentzian unimodular and non-unimodular Lie groups. In section 4, we recall the definition of the Bott connection $\nabla^{B_3}$ with the first distribution, and then give left-invariant Ricci collineations associated to the Bott connection $\nabla^{B_3}$ on three-dimensional Lorentzian unimodular and non-unimodular Lie groups.


\vskip 1 true cm

\section{ Invariant Ricci collineations associated to the Bott connection on three-dimensional Lorentzian Unimodular Lie groups with the first distribution}

\indent Calvaruso and Cordero classified three-dimensional Lorentzian Lie groups in \cite{Ca1,CP}(see Theorem 2.1 and Theorem 2.2 in \cite{W}). Throughout this paper, we shall denote the connected by $\{G_i\}_{i=1,\cdots,7}$, simply connected three-dimensional Lie group equipped with a left-invariant Lorentzian metric $g$ and having Lie algebra $\{\mathfrak{g}\}_{i=1,\cdots,7.}$.\\
\indent Then we recall the definition of the Bott connection $\nabla^{B_1}$. Let $M$ be a smooth manifold, and let $TM=span\{e_1,e_2,e_3\}$, then took the frst distribution: $F_1=span\{e_1,e_2\}$ and ${F_1}^\bot=span\{e_3\}$, where ${e_1,e_2,e_3}$ is a pseudo-orthonormal basis, with $e_3$ timelike.
The  Bott connection $\nabla^{B_1}$ is defined as follows: (see \cite{F}, \cite{J}, \cite{RK})
\begin{eqnarray}
\nabla^{B_1}_XY=
       \begin{cases}
        \pi_{F_1}(\nabla^L_XY),~~~&X,Y\in\Gamma^\infty(F_1) \\[2pt]
       \pi_{F_1}([X,Y]),~~~&X\in\Gamma^\infty({F_1}^\bot),Y\in\Gamma^\infty(F_1)\\[2pt]
       \pi_{{F_1}^\bot}([X,Y]),~~~&X\in\Gamma^\infty(F_1),Y\in\Gamma^\infty({F_1}^\bot)\\[2pt]
       \pi_{{F_1}^\bot}(\nabla^L_XY),~~~&X,Y\in\Gamma^\infty({F_1}^\bot)\\[2pt]
       \end{cases}
\end{eqnarray}
where $\pi_{F_1}$ and $\pi_{F_1}^\bot$ are respectively the projection on $F_1$ and ${F_1}^\bot$, $\nabla^L$ is the Levi-Civita connection of $G_i$ .\\
We define the curvature tensor of the Bott connection $\nabla^{B_1}$
\begin{equation}
R^{B_1}(X,Y)Z=\nabla^{B_1}_X\nabla^{B_1}_YZ-\nabla^{B_1}_Y\nabla^{B_1}_XZ-\nabla^{B_1}_{[X,Y]}Z.
\end{equation}
The Ricci tensor of $(G_i,g)$ associated to the Bott connection $\nabla^{B_1}$ is defined by
\begin{equation}
\rho^{B_1}(X,Y)=-g(R^{B_1}(X,e_1)Y,e_1)-g(R^{B_1}(X,e_2)Y,e_2)+g(R^{B_1}(X,e_3)Y,e_3).
\end{equation}
Let
\begin{equation}
Ric^{B_1}(X,Y)=\frac{\rho^{B_1}(X,Y)+\rho^{B_1}(Y,X)}{2}.
\end{equation}
We define:
\begin{equation}
(L_VRic^{B_1})(X,Y):=V[Ric^{B_1}(X,Y)]-Ric^{B_1}([V,X],Y)-Ric^{B_1}(X,[V,Y])
\end{equation}
for vector $X,Y,V$.
\begin{thm}$(G_i,g)$ admits left-invariant Ricci collineations  associated to the Bott connection $\nabla^{B_1}$ if and only if it satisfies
\begin{equation}
(L_VRic^{B_1})(X,Y)=0,
\end{equation}
where $V=\lambda_1e_1+\lambda_2e_2+\lambda_3e_3$ is a left-invariant vector field and $\lambda_1,\lambda_2,\lambda_3$ are real numbers.
\end{thm}

\vskip 0.5 true cm
\noindent{\bf 2.1 Invariant Ricci collineations of $G_1$ associated to the Bott connection $\nabla^{B_1}$}\\
\vskip 0.5 true cm
 By \cite{W}, we have the following Lie algebra of $G_1$ satisfies
\begin{equation}
[e_1,e_2]=\alpha e_1-\beta e_3,~~[e_1,e_3]=-\alpha e_1-\beta e_2,~~[e_2,e_3]=\beta e_1+\alpha e_2+\alpha e_3,~~\alpha\neq 0.
\end{equation}
where $e_1,e_2,e_3$ is a pseudo-orthonormal basis, with $e_3$ timelike.
\begin{lem}(\cite{Wu}) The Ricci tensor of  $(G_1,g)$ associated to the Bott connection $\nabla^{B_1}$ is determined by
\begin{align}
&Ric^{B_1}(e_1,e_1)=-(\alpha^2+\beta^2),~~~Ric^{B_1}(e_1,e_2)=\alpha\beta,~~~Ric^{B_1}(e_1,e_3)=-\frac{\alpha\beta}{2},\\\notag
&Ric^{B_1}(e_2,e_2)=-(\alpha^2+\beta^2),~~~Ric^{B_1}(e_2,e_3)=\frac{\alpha^2}{2},~~~Ric^{B_3}(e_3,e_3)=0.
\end{align}
\end{lem}
By (2.5) and Lemma 2.2, we have
\begin{lem}
\begin{align}
&(L_VRic^{B_1})(e_1,e_1)=-\alpha(2\alpha^2+\beta^2)\lambda_2+2\alpha^3\lambda_3,\\\notag
&(L_VRic^{B_1})(e_1,e_2)=\alpha(\alpha^2+\dfrac{\beta^2}{2})\lambda_1+\dfrac{\alpha^2\beta}{2}\lambda_2-\dfrac{\alpha^2\beta}{2}\lambda_3,\\\notag
&(L_VRic^{B_1})(e_1,e_3)=-\alpha^3\lambda_1+\beta^3\lambda_2,\\\notag
&(L_VRic^{B_1})(e_2,e_2)=-\alpha^2\beta\lambda_1-\alpha^3\lambda_3,\\\notag
&(L_VRic^{B_1})(e_2,e_3)=\beta(\dfrac{\alpha^2}{2}-\beta^2)\lambda_1+\dfrac{\alpha^3}{2}\lambda_2+\dfrac{\alpha}{2}(\alpha^2-\beta^2)\lambda_3,\\\notag
&(L_VRic^{B_1})(e_3,e_3)=\alpha(\beta^2-\alpha^2)\lambda_2.
\end{align}
\end{lem}
Then, if a left-invariant vector field $V$ is a  Ricci collineation associated to the Bott connection $\nabla^{B_1}$, by Lemma 2.3 and Theorem 2.1, we have the following equations:
\begin{eqnarray}
       \begin{cases}
       -\alpha(2\alpha^2+\beta^2)\lambda_2+2\alpha^3\lambda_3=0\\[2pt]
       \alpha(\alpha^2+\dfrac{\beta^2}{2})\lambda_1+\dfrac{\alpha^2\beta}{2}\lambda_2-\dfrac{\alpha^2\beta}{2}\lambda_3=0\\[2pt]
       -\alpha^3\lambda_1+\beta^3\lambda_2=0\\[2pt]
       -\alpha^2\beta\lambda_1-\alpha^3\lambda_3=0\\[2pt]
       \beta(\dfrac{\alpha^2}{2}-\beta^2)\lambda_1+\dfrac{\alpha^3}{2}\lambda_2+\dfrac{\alpha}{2}(\alpha^2-\beta^2)\lambda_3=0\\[2pt]
       \alpha(\beta^2-\alpha^2)\lambda_2=0\\[2pt]
       \end{cases}
\end{eqnarray}
By solving (2.10) , we get
\begin{thm}
$(G_1, g, V)$ does not admit left-invariant Ricci collineations associated to the Bott connection $\nabla^{B_1}$.
\end{thm}
\begin{Proof}
We know that $\alpha\neq 0$. By the sixth equation, we get $(\beta^2-\alpha^2)\lambda_2=0$, so\\
case {\rm 1)} If $\alpha=\beta$, by {\rm (2.10)},
\begin{eqnarray}
       \begin{cases}
       -3\lambda_2+2\lambda_3=0\\[2pt]
       3\lambda_1+\lambda_2-\lambda_3=0\\[2pt]
       -\lambda_1+\lambda_2=0\\[2pt]
       \lambda_1+\lambda_3=0\\[2pt]
       \end{cases}
\end{eqnarray}
then we get $\lambda_1=\lambda_2=\lambda_3=0$.\\
case {\rm 2)} If $\alpha=-\beta$, by {\rm (2.10)},
\begin{eqnarray}
       \begin{cases}
       -3\lambda_2+2\lambda_3=0\\[2pt]
       3\lambda_1-\lambda_2+\lambda_3=0\\[2pt]
       \lambda_1+\lambda_2=0\\[2pt]
       \lambda_1-\lambda_3=0\\[2pt]
       \end{cases}
\end{eqnarray}
then we get $\lambda_1=\lambda_2=\lambda_3=0$.\\
case {\rm 3)} If $\lambda_2=0$, by the third equation, we get $\alpha^3\lambda_1=0$, i.e. $\lambda_1=0$, then by the fourth equation, we get $\alpha^3\lambda_3=0$, i.e. $\lambda_3=0$, so $\lambda_1=\lambda_2=\lambda_3=0$.\\
So we have no non-trivial solution.
\end{Proof}

\vskip 0.5 true cm
\noindent{\bf 2.2 Invariant Ricci collineations of $G_2$ associated to the Bott connection $\nabla^{B_1}$}\\
\vskip 0.5 true cm
 By \cite{W}, we have the following Lie algebra of $G_2$ satisfies
\begin{equation}
[e_1,e_2]=\gamma e_2-\beta e_3,~~[e_1,e_3]=-\beta e_2-\gamma e_3,~~[e_2,e_3]=\alpha e_1,~~\gamma\neq 0.
\end{equation}
where $e_1,e_2,e_3$ is a pseudo-orthonormal basis, with $e_3$ timelike.
\begin{lem}(\cite{Wu}) The Ricci tensor of  $(G_2,g)$ associated to the Bott connection $\nabla^{B_1}$ is determined by
\begin{align}
&Ric^{B_1}(e_1,e_1)=-(\beta^2+\gamma^2),~~~Ric^{B_1}(e_1,e_2)=0,~~~Ric^{B_1}(e_1,e_3)=0,\\\notag
&Ric^{B_1}(e_2,e_2)=-(\gamma^2+\alpha\beta),~~~Ric^{B_1}(e_2,e_3)=-\frac{\alpha\gamma}{2},~~~Ric^{B_3}(e_3,e_3)=0.
\end{align}
\end{lem}
By (2.5) and Lemma 2.5, we have
\begin{lem}
\begin{align}
&(L_VRic^{B_1})(e_1,e_1)=0,\\\notag
&(L_VRic^{B_1})(e_1,e_2)=-\gamma(\gamma^2+\dfrac{\alpha\beta}{2})\lambda_2+\gamma^2(\beta-\dfrac{\alpha}{2})\lambda_3,\\\notag
&(L_VRic^{B_1})(e_1,e_3)=\alpha(\beta^2+\dfrac{\gamma^2}{2})\lambda_2+\dfrac{\alpha\beta\gamma}{2}\lambda_3,\\\notag
&(L_VRic^{B_1})(e_2,e_2)=\gamma(2\gamma^2+\alpha\beta)\lambda_1,\\\notag
&(L_VRic^{B_1})(e_2,e_3)=-\beta(\gamma^2+\alpha\beta)\lambda_1,\\\notag
&(L_VRic^{B_1})(e_3,e_3)=-\alpha\beta\gamma\lambda_1.
\end{align}
\end{lem}
\indent Then, if a left-invariant vector field $V$ is a  Ricci collineation associated to the Bott connection $\nabla^{B_1}$, by Lemma 2.4 and Theorem 2.1, we have the following equations:
\begin{eqnarray}
       \begin{cases}
       -\gamma(\gamma^2+\dfrac{\alpha\beta}{2})\lambda_2+\gamma^2(\beta-\dfrac{\alpha}{2})\lambda_3=0\\[2pt]
       \alpha(\beta^2+\dfrac{\gamma^2}{2})\lambda_2+\dfrac{\alpha\beta\gamma}{2}\lambda_3=0\\[2pt]
       \gamma(2\gamma^2+\alpha\beta)\lambda_1=0\\[2pt]
       -\beta(\gamma^2+\alpha\beta)\lambda_1=0\\[2pt]
       -\alpha\beta\gamma\lambda_1=0\\[2pt]
       \end{cases}
\end{eqnarray}
\indent By solving (2.16) , we get
\begin{thm}
$(G_2, g, V)$ admits left-invariant Ricci collineations associated to the Bott connection $\nabla^{B_1}$ if and only if one of the following holds:\\
{\rm (1)} $\alpha=0,\beta= 0,\gamma\neq 0$,\\
{\rm (2)} $\alpha=0,\beta\neq 0,\gamma\neq 0$,\\
{\rm (3)} $\alpha\neq 0,\beta\neq 0,\gamma\neq 0, 2\beta^3-2\alpha\beta^2-\dfrac{\alpha\gamma^2}{2}=0.$\\
\indent Moreover, in these cases, we have\\
{\rm (1)} $\mathscr{V}_{\mathscr{R}{C}}=<e_2>$.\\
{\rm (2)} $\mathscr{V}_{\mathscr{R}{C}}=<\dfrac{\beta}{\gamma}e_2+e_3>$.\\
{\rm (3)} $\mathscr{V}_{\mathscr{R}{C}}=<-\dfrac{\beta\gamma}{2\beta^2+\gamma^2}e_2+e_3>$.\\
where $\mathscr{V}_{\mathscr{R}{C}}$ is the vector space of left-invariant Ricci collineations on $(G_2, g, V)$.
\end{thm}
\begin{Proof}
We know that $\gamma\neq 0$, by the fifth equation, we have $\alpha\beta\lambda_1=0$. \\
case {\rm 1)} If $\lambda_1\neq 0$, then $\alpha\beta=0$. Note that by the third equation, we have
\begin{eqnarray}
\gamma(2\gamma^2+\alpha\beta)\lambda_1=2\gamma^3\lambda_1= 0,\notag
\end{eqnarray}
i.e. $\gamma=0$. This is an contracdiction.\\
case {\rm 2)} If $\lambda_1=0$, then the third,foyrth,fifth equation trivially holds.\\
case {\rm 2-1)} If $\alpha=0$, then by the first equation, we have $-\gamma^3\lambda_2+\gamma^2\beta\lambda_3=0$, i.e.
\begin{eqnarray}
-\gamma\lambda_2+\beta\lambda_3=0.\notag
 \end{eqnarray}
case {\rm 2-1-1)} If $\beta=0$, then $\lambda_2= 0$. We get {\rm (1)}.\\
case {\rm 2-1-2)} If $\beta\neq 0$, then $-\gamma\lambda_2+\beta\lambda_3=0$, i.e. $\lambda_2=\dfrac{\beta}{\gamma}\lambda_3$. We get {\rm (2)}.\\
case {\rm 2-2)} If $\alpha\neq 0,$ by {\rm (2.16)}, we have
\begin{eqnarray}
\begin{cases}
(\gamma^2+\alpha\beta)\lambda_2+\gamma(\beta-\dfrac{\alpha}{2})\lambda_3=0\\[2pt]
(\beta^2+\dfrac{\gamma^2}{2})\lambda_2+\dfrac{\beta\gamma}{2}\lambda_3=0\\[2pt]
\end{cases}
\end{eqnarray}
case {\rm 2-2-1)} If $\beta=0$,  by {\rm (2.17)}, we have
\begin{eqnarray}
\begin{cases}
\gamma^2\lambda_2+\dfrac{\alpha\gamma}{2}\lambda_3=0\\[2pt]
\dfrac{\alpha\gamma^2}{2}\lambda_2=0\\[2pt]
\end{cases}
\end{eqnarray}
then we have $\lambda_2=0,\lambda_3=0,$ i.e. $\lambda_1=\lambda_2=\lambda_3=0$.\\
case {\rm 2-2-2)} If $\beta\neq 0$, solve {\rm (2.16)} we have
\begin{eqnarray}
2\beta^3-2\alpha\beta^2-\dfrac{\alpha\gamma^2}{2}=0,\notag
\end{eqnarray}
then $\lambda_2=
-\dfrac{\beta\gamma}{2\beta^2+\gamma^2}\lambda_3$. We get {\rm (3)}\\
case {\rm 2-2)} If $\alpha=0$, by {\rm (2.16)}, we have
\begin{eqnarray}
-\gamma\lambda_2+\beta\lambda_3=0.\notag
\end{eqnarray}

\end{Proof}

\vskip 0.5 true cm
\noindent{\bf 2.3 Invariant Ricci collineations of $G_3$ associated to the Bott connection $\nabla^{B_1}$}\\
\vskip 0.5 true cm
 By \cite{W}, we have the following Lie algebra of $G_3$ associated to the Bott connection $\nabla^{B_1}$ satisfies
\begin{equation}
[e_1,e_2]=-\gamma e_3,~~[e_1,e_3]=-\beta e_2,~~[e_2,e_3]=\alpha e_1.
\end{equation}
where $e_1,e_2,e_3$ is a pseudo-orthonormal basis, with $e_3$ timelike.
\begin{lem}(\cite{Wu}) The Ricci tensor of  $(G_2,g)$ is determined by
\begin{align}
&Ric^{B_1}(e_1,e_1)=-\beta\gamma,~~~Ric^{B_1}(e_1,e_2)=0,~~~Ric^{B_1}(e_1,e_3)=0,\\\notag
&Ric^{B_1}(e_2,e_2)=-\alpha\gamma,~~~Ric^{B_1}(e_2,e_3)=0,~~~Ric^{B_3}(e_3,e_3)=0.
\end{align}
\end{lem}
By (2.5) and Lemma 2.8, we have
\begin{lem}
\begin{align}
&(L_VRic^{B_1})(e_1,e_1)=0,~~~(L_VRic^{B_1})(e_1,e_2)=0,\\\notag
&(L_VRic^{B_1})(e_1,e_3)=\alpha\beta\gamma\lambda_2,~~~(L_VRic^{B_1})(e_2,e_2)=0,\\\notag
&(L_VRic^{B_1})(e_2,e_3)=-\alpha\beta\gamma\lambda_1,~~~(L_VRic^{B_1})(e_3,e_3)=0.
\end{align}
\end{lem}
\indent Then, if a left-invariant vector field $V$ is a  Ricci collineation associated to the Bott connection $\nabla^{B_1}$, by Lemma 2.7 and Theorem 2.1, we have the following equations:
\begin{eqnarray}
       \begin{cases}
       \alpha\beta\gamma\lambda_2=0\\[2pt]
       -\alpha\beta\gamma\lambda_1=0\\[2pt]
       \end{cases}
\end{eqnarray}
\indent By solving (2.22) , we get
\begin{thm}
$(G_3, g, V)$ admits left-invariant Ricci collineations associated to the Bott connection $\nabla^{B_1}$ if and only if one of the following holds:\\
{\rm (1)} $\alpha\beta\gamma=0$ ,\\
{\rm (2)} $\alpha\neq 0,\beta\neq 0,\gamma\neq 0$.\\
\indent Moreover, in these cases, we have\\
{\rm (1)} $\mathscr{V}_{\mathscr{R}{C}}=<e_1,e_2,e_3>$,\\
{\rm (2)} $\mathscr{V}_{\mathscr{R}{C}}=<e_3>$.\\
where $\mathscr{V}_{\mathscr{R}{C}}$ is the vector space of left-invariant Ricci collineations on $(G_3, g, V)$.
\end{thm}
\begin{Proof}
case {\rm 1)} If $\alpha\beta\gamma=0$, {\rm (2.22)} trivially holds. We get {\rm (1)}.\\
case {\rm 2)} If $\alpha\beta\gamma\neq 0$, i.e. $\alpha\neq 0,\beta\neq 0,\gamma\neq 0$, then $\lambda_1=\lambda_2=0$. We get {\rm (2)}.\\
\end{Proof}

\vskip 0.5 true cm
\noindent{\bf 2.4 Invariant Ricci collineations of $G_4$ associated to the Bott connection $\nabla^{B_1}$}\\
\vskip 0.5 true cm
 By \cite{W}, we have the following Lie algebra of $G_4$ satisfies
\begin{equation}
[e_1,e_2]=-e_2+(2\eta-\beta)e_3,~~\eta=1~{\rm or}-1,~~[e_1,e_3]=-\beta e_2+ e_3,~~[e_2,e_3]=\alpha e_1.
\end{equation}
where $e_1,e_2,e_3$ is a pseudo-orthonormal basis, with $e_3$ timelike.
\begin{lem}(\cite{Wu}) The Ricci tensor of  $(G_4,g)$ associated to the Bott connection $\nabla^{B_1}$ is determined by
\begin{align}
&Ric^{B_1}(e_1,e_1)=-(\beta-\eta)^2,~~~Ric^{B_1}(e_1,e_2)=0,~~~Ric^{B_1}(e_1,e_3)=0,\\\notag
&Ric^{B_1}(e_2,e_2)=2\alpha\eta-\alpha\beta-1,~~~Ric^{B_1}(e_2,e_3)=\dfrac{\alpha}{2},~~~Ric^{B_3}(e_3,e_3)=0.
\end{align}
\end{lem}
By (2.5) and Lemma(2.11), we have
\begin{lem}
\begin{align}
&(L_VRic^{B_1})(e_1,e_1)=0,\\\notag
&(L_VRic^{B_1})(e_1,e_2)=(-\alpha\eta+\dfrac{\alpha\beta}{2}+1)\lambda_2+(\beta+\dfrac{\alpha}{2}-\alpha\eta^2)\lambda_3,\\\notag
&(L_VRic^{B_1})(e_1,e_3)=\alpha[(\beta-\eta)^2-\dfrac{1}{2}]\lambda_2-\dfrac{\alpha\beta}{2}\lambda_3,\\\notag
&(L_VRic^{B_1})(e_2,e_2)=(2\alpha\eta-\alpha\beta-2)\lambda_1,\\\notag
&(L_VRic^{B_1})(e_2,e_3)=\beta(2\alpha\eta-\alpha\beta-1)\lambda_1,\\\notag
&(L_VRic^{B_1})(e_3,e_3)=\alpha\beta\lambda_1.
\end{align}
\end{lem}
\indent Then, if a left-invariant vector field $V$ is a  Ricci collineation associated to the Bott connection $\nabla^{B_1}$, by Lemma 2.10 and Theorem 2.1, we have the following equations:
\begin{eqnarray}
       \begin{cases}
       (-\alpha\eta+\dfrac{\alpha\beta}{2}+1)\lambda_2+(\beta+\dfrac{\alpha}{2}-\alpha\eta^2)\lambda_3=0\\[2pt]
       \alpha[(\beta-\eta)^2-\dfrac{1}{2}]\lambda_2-\dfrac{\alpha\beta}{2}\lambda_3=0\\[2pt]
       (2\alpha\eta-\alpha\beta-2)\lambda_1=0\\[2pt]
       \beta(2\alpha\eta-\alpha\beta-1)\lambda_1=0\\[2pt]
       \alpha\beta\lambda_1=0\\[2pt]
       \end{cases}
\end{eqnarray}
\indent By solving (2.26) , we get
\begin{thm}
$(G_4, g, V)$ admits left-invariant Ricci collineations associated to the Bott connection $\nabla^{B_1}$ if and only if one of the following holds:\\
{\rm (1)} $\alpha\neq 0,\beta=0,\eta=1~{\rm or}-1,\alpha\eta=1, $\\
{\rm (2)} $\alpha=0,\beta=0,\eta=1~{\rm or}-1,$\\
{\rm (3)} $\alpha= 0,\beta\neq 0,\eta=1~{\rm or}-1,$\\
{\rm (4)} $\alpha\neq 0,\beta\neq 0,\eta=1~{\rm or}-1,\alpha-4\beta=0,$\\
{\rm (5)} $\alpha\neq 0,\beta\neq 0,\eta=1~{\rm or}-1,\beta=\eta.$\\
\indent Moreover, in these cases, we have\\
{\rm (1)} $\mathscr{V}_{\mathscr{R}{C}}=<e_1>$,\\
{\rm (2)} $\mathscr{V}_{\mathscr{R}{C}}=<e_3>$,\\
{\rm (3)} $\mathscr{V}_{\mathscr{R}{C}}=<-\beta e_2+e_3>$,\\
{\rm (4)} $\mathscr{V}_{\mathscr{R}{C}}=<e_2+(2\beta+\dfrac{1}{\beta}-4\eta)e_3>$,\\
{\rm (5)} $\mathscr{V}_{\mathscr{R}{C}}=<-\eta e_2+e_3>$.\\
where $\mathscr{V}_{\mathscr{R}{C}}$ is the vector space of left-invariant Ricci collineations on $(G_4, g, V)$.
\end{thm}
\begin{Proof}
We know that $\eta=1~{\rm or}-1$.\\
case {\rm 1)} If $\lambda_1\neq 0$ by the sixth equation, $\alpha\beta=0$.\\
case {\rm 1-1)} If $\alpha=0$, by the third equation we have $-2\lambda_1=0$, i.e. $\lambda_1=0$. This is a contradiction.\\
case {\rm 1-2)} If $\alpha\neq 0$, then $\beta=0,$ by the third equation we have $2\alpha\eta-2=0$, i.e. $\alpha\eta=1$.
By {\rm (2.26)}, we have
\begin{eqnarray}
\begin{cases}
-\dfrac{\alpha}{2}\lambda_3=0\\[2pt]
\dfrac{\alpha}{2}\lambda_2=0\\[2pt]
\end{cases}
\end{eqnarray}
i.e. $\lambda_2=\lambda_3=0$. We get {\rm (1)}.\\
case {\rm 2)} If $\lambda_1= 0,$ by {\rm (2.26)},
\begin{eqnarray}
\begin{cases}
(-\alpha\eta+\dfrac{\alpha\beta}{2}+1)\lambda_2+(\beta-\dfrac{\alpha}{2})\lambda_3=0\\[2pt]
\alpha[(\beta-\eta)^2-\dfrac{1}{2}]\lambda_2-\dfrac{\alpha\beta}{2}\lambda_3=0\\[2pt]
\end{cases}
\end{eqnarray}
case {\rm 2-1)} If $\alpha=0$, then $\lambda_2+\beta\lambda_3=0$.\\
case {\rm 2-1-1)} If $\beta= 0$, then $\lambda_2=0$. We get {\rm (2)}.\\
case {\rm 2-1-2)} If $\beta\neq 0$, then $\lambda_2=-\beta\lambda_3$. We get {\rm (3)}.\\
case {\rm 2-2)} If $\alpha\neq 0,$ we get
\begin{eqnarray}
\begin{cases}
(-\alpha\eta+\dfrac{\alpha\beta}{2}+1)\lambda_2+(\beta-\dfrac{\alpha}{2})\lambda_3=0\\[2pt]
[(\beta-\eta)^2-\dfrac{1}{2}]\lambda_2-\dfrac{\beta}{2}\lambda_3=0\\[2pt]
\end{cases}
\end{eqnarray}
case {\rm 2-2-1)} If $\beta=0$, we have
\begin{eqnarray}
\begin{cases}
(-\alpha\eta+1)\lambda_2-\dfrac{\alpha}{2}\lambda_3=0\\[2pt]
\dfrac{1}{2}\lambda_2=0\\[2pt]
\end{cases}
\end{eqnarray}
then $\lambda_2=0, \lambda_3=0$, i.e. $\lambda_1=\lambda_2=\lambda_3=0$.\\
case {\rm 2-2-2)} If $\beta\neq 0$, we have $[[(\beta-\eta)^2-\dfrac{1}{2}](\beta-\dfrac{\alpha}{2})+\dfrac{\beta}{2}(-\alpha\eta+\dfrac{\alpha\beta}{2}+1)]\lambda_2=0$, i.e.
\begin{eqnarray}
(4\beta-\alpha)(\beta-\eta)^2\lambda_2=0.\notag
\end{eqnarray}
If $\lambda_2=0$, then $\lambda_3=0$, i.e. $\lambda_1=\lambda_2=\lambda_3=0$.\\
If $\alpha=4\beta$, then $\lambda_3=(2\beta+\dfrac{1}{\beta}-4\eta)\lambda_2$. We get {\rm (4)}.\\
If $\beta=\eta$, then $\lambda_2=-\beta\lambda_3$. We get {\rm (5)}.\\
\end{Proof}

\vskip 0.5 true cm
\noindent{\bf 2.5 Invariant Ricci collineations of $G_5$ associated to the Bott connection $\nabla^{B_1}$}\\
\vskip 0.5 true cm
 By \cite{W}, we have the following Lie algebra of $G_5$ satisfies
\begin{equation}
[e_1,e_2]=0,~~[e_1,e_3]=\alpha e_1+\beta e_2,~~[e_2,e_3]=\gamma e_1+\delta e_2,~~\alpha+\delta\neq 0,~~\alpha\gamma+\beta\delta=0.
\end{equation}
where $e_1,e_2,e_3$ is a pseudo-orthonormal basis, with $e_3$ timelike.
\begin{lem}(\cite{Wu}) The Ricci tensor of  $(G_5,g)$ associated to the Bott connection $\nabla^{B_1}$ is determined by
\begin{align}
&Ric^{B_1}(e_i,e_j)=0.
\end{align}
for any pairs (i,j).
\end{lem}
By (2.5) and Lemma(2.14), we have
\begin{lem}
\begin{align}
&(L_VRic^{B_1})(e_i,e_j)=0.
\end{align}
\end{lem}
\indent Then we have
\begin{thm}
Any left-invariant vector field on $(G_5, g, V)$ is a left-invariant Ricci collineations associated to the Bott connection $\nabla^{B_1}$ .\\
\end{thm}

\vskip 0.5 true cm
\noindent{\bf 2.6 Invariant Ricci collineations of $G_6$ associated to the Bott connection $\nabla^{B_1}$}\\
\vskip 0.5 true cm
 By \cite{W}, we have the following Lie algebra of $G_6$ satisfies
\begin{equation}
[e_1,e_2]=\alpha e_2+\beta e_3,~~[e_1,e_3]=\gamma e_2+\delta e_3,~~[e_2,e_3]=0,~~\alpha+\delta\neq 0£¬~~\alpha\gamma-\beta\delta=0.
\end{equation}
where $e_1,e_2,e_3$ is a pseudo-orthonormal basis, with $e_3$ timelike.
\begin{lem}(\cite{Wu}) The Ricci tensor of  $(G_6,g)$ associated to the Bott connection $\nabla^{B_1}$ is determined by
\begin{align}
&Ric^{B_1}(e_1,e_1)=-(\alpha^2+\beta\gamma),~~~Ric^{B_1}(e_1,e_2)=0,~~~Ric^{B_1}(e_1,e_3)=0,\\\notag
&Ric^{B_1}(e_2,e_2)=-\alpha^2,~~~Ric^{B_1}(e_2,e_3)=0,~~~Ric^{B_3}(e_3,e_3)=0.
\end{align}
\end{lem}
By (2.5) and Lemma 2.17, we have
\begin{lem}
\begin{align}
&(L_VRic^{B_1})(e_1,e_1)=0,~~~(L_VRic^{B_1})(e_1,e_2)=-\alpha^3\lambda_2-\alpha^2\gamma\lambda_3,\\\notag
&(L_VRic^{B_1})(e_1,e_3)=0~~~(L_VRic^{B_1})(e_2,e_2)=2\alpha^3\lambda_1,\\\notag
&(L_VRic^{B_1})(e_2,e_3)=\alpha^2\gamma\lambda_1,~~~(L_VRic^{B_1})(e_3,e_3)=0.
\end{align}
\end{lem}
\indent Then, if a left-invariant vector field $V$ is a  Ricci collineation associated to the Bott connection $\nabla^{B_1}$, by Lemma 2.17 and Theorem 2.1, we have the following equations:
\begin{eqnarray}
       \begin{cases}
       -\alpha^3\lambda_2-\alpha^2\gamma\lambda_3=0\\[2pt]
       2\alpha^3\lambda_1=0\\[2pt]
       \alpha^2\gamma\lambda_1=0\\[2pt]
       \end{cases}
\end{eqnarray}
\indent By solving (2.37) , we get
\begin{thm}
$(G_6, g, V)$ admits left-invariant Ricci collineations associated to the Bott connection $\nabla^{B_1}$ if and only if one of the following holds:\\
{\rm (1)} $\alpha= 0,\beta=0,\delta\neq 0, $\\
{\rm (2)} $\alpha\neq 0,\gamma=0,\alpha+\delta\neq 0,\beta\delta= 0,$\\
{\rm (3)} $\alpha\neq 0,\gamma\neq 0,\alpha+\delta\neq 0,\alpha\gamma-\beta\delta= 0.$\\
\indent Moreover, in these cases, we have\\
{\rm (1)} $\mathscr{V}_{\mathscr{R}{C}}=<e_1,e_2,e_3>$,\\
{\rm (2)} $\mathscr{V}_{\mathscr{R}{C}}=<e_3>$,\\
{\rm (3)} $\mathscr{V}_{\mathscr{R}{C}}=<-\dfrac{\gamma}{\alpha}e_2+e_3>$.\\
where $\mathscr{V}_{\mathscr{R}{C}}$ is the vector space of left-invariant Ricci collineations on $(G_6, g, V)$.
\end{thm}

\begin{Proof}
We know that $\alpha+\delta\neq 0,\alpha\gamma-\beta\delta= 0$.\\
case {\rm 1)} If $\alpha=0$, then $\delta\neq 0, \beta=0$. {\rm (2.37)} trivially holds.\\
case {\rm 2)} If $\alpha\neq 0,$ by the second equation we have $\lambda_1=0$. By {\rm (2.37)}, we have
\begin{eqnarray}
 \alpha\lambda_2+\gamma\lambda_3=0.\notag
 \end{eqnarray}
case {\rm 2-1)} If $\gamma=0$, then $\lambda_2=0,\beta\delta= 0$.We get {\rm (2)}.\\
case {\rm 2-2)} If $\gamma\neq 0$, then $\lambda_2=-\dfrac{\gamma}{\alpha}\lambda_3$. We get {\rm (3)}.\\
\end{Proof}

\vskip 0.5 true cm
\noindent{\bf 2.7 Invariant Ricci collineations of $G_7$ associated to the Bott connection $\nabla^{B_1}$}\\
\vskip 0.5 true cm
 By \cite{W}, we have the following Lie algebra of $G_7$ satisfies
\begin{equation}
[e_1,e_2]=-\alpha e_1-\beta e_2-\beta e_3,~~[e_1,e_3]=\alpha e_1+\beta e_2+\beta e_3,~~[e_2,e_3]=\gamma e_1+\delta e_2+\delta e_3,,~~\alpha+\delta\neq 0,~~\alpha\gamma=0.
\end{equation}
where $e_1,e_2,e_3$ is a pseudo-orthonormal basis, with $e_3$ timelike.
\begin{lem}(\cite{Wu}) The Ricci tensor of  $(G_7,g)$ associated to the Bott connection $\nabla^{B_1}$ is determined by
\begin{align}
&Ric^{B_1}(e_1,e_1)=-\alpha^2,~~~Ric^{B_1}(e_1,e_2)=\dfrac{\beta(\delta-\alpha)}{2},\\\notag
&Ric^{B_1}(e_1,e_3)=\beta(\alpha+\delta),Ric^{B_1}(e_2,e_2)=-(\alpha^2+\beta^2+\beta\gamma),\\\notag
&Ric^{B_1}(e_2,e_3)=\delta^2+\dfrac{\beta\gamma+\alpha\delta}{2},~~~Ric^{B_3}(e_3,e_3)=0.
\end{align}
\end{lem}
By (2.5) and Lemma 2.20, we have
\begin{lem}
\begin{align}
&(L_VRic^{B_1})(e_1,e_1)=(2\alpha^3-3\beta^2\delta-\alpha\beta^2)\lambda_2-(2\alpha^3-3\beta^2\delta-\alpha\beta^2)\lambda_3,\\\notag
&(L_VRic^{B_1})(e_1,e_2)=(-\alpha^3+\dfrac{3\beta^2\delta}{2}-\dfrac{\alpha\beta^2}{2})\lambda_1+\beta(\dfrac{3\alpha^2}{2}+\beta^2-\delta^2-\alpha\delta+\dfrac{\beta\gamma}{2})\lambda_2\\\notag
&\hspace{3.5cm}+\beta(-\dfrac{3\alpha^2}{2}-\beta^2+\dfrac{5\delta^2}{2}+\dfrac{3\alpha\delta}{2}-\dfrac{\beta\gamma}{2})\lambda_3,\\\notag
&(L_VRic^{B_1})(e_1,e_3)=(\alpha^3-\dfrac{3\beta^2\delta}{2}+\dfrac{\alpha\beta^2}{2})\lambda_1-\beta(\alpha^2+\dfrac{5\delta^2}{2}+2\alpha\delta+\dfrac{\beta\gamma}{2})\lambda_2\\\notag
&\hspace{3.5cm}+\beta(\alpha^2+\delta^2+\dfrac{3\alpha\delta}{2}+\dfrac{\beta\gamma}{2})\lambda_3,\\\notag
&(L_VRic^{B_1})(e_2,e_2)=\beta(-3\alpha^2-2\beta^2+2\delta^2+2\alpha\delta-\beta\gamma)\lambda_1+\delta(-2\alpha^2-2\beta^2+2\delta^2+\alpha\delta)\lambda_3,\\\notag
&(L_VRic^{B_1})(e_2,e_3)=\beta(\dfrac{5\alpha^2}{2}+\beta^2+\dfrac{\alpha\delta}{2}+\beta\gamma)\lambda_1+\delta(\alpha^2+\beta^2-\delta^2-\dfrac{\alpha\delta}{2})\lambda_2,\\\notag
&\hspace{3.5cm}+\delta(\delta^2+\dfrac{\alpha\delta}{2}+\dfrac{3\beta\gamma}{2})\lambda_3,\\\notag
&(L_VRic^{B_1})(e_3,e_3)=-\beta(2\alpha^2+2\delta^2+3\alpha\delta+\beta\gamma)\lambda_1-\delta(-2\delta^2+\alpha\delta+3\beta\delta)\lambda_3.
\end{align}
\end{lem}
\indent Then, if a left-invariant vector field $V$ is a  Ricci collineation associated to the Bott connection $\nabla^{B_1}$, by Lemma 2.19 and Theorem 2.1, we have the following equations:
\begin{eqnarray}
       \begin{cases}
       (2\alpha^3-3\beta^2\delta-\alpha\beta^2)\lambda_2-(2\alpha^3-3\beta^2\delta-\alpha\beta^2)\lambda_3=0\\[2pt]
       (-\alpha^3+\dfrac{3\beta^2\delta}{2}-\dfrac{\alpha\beta^2}{2})\lambda_1+\beta(\dfrac{3\alpha^2}{2}+\beta^2-\delta^2-\alpha\delta+\dfrac{\beta\gamma}{2})\lambda_2\\[2pt]
       \hspace{1cm}+\beta(-\dfrac{3\alpha^2}{2}-\beta^2+\dfrac{5\delta^2}{2}+\dfrac{3\alpha\delta}{2}-\dfrac{\beta\gamma}{2})\lambda_3=0\\[2pt]
       (\alpha^3-\dfrac{3\beta^2\delta}{2}+\dfrac{\alpha\beta^2}{2})\lambda_1-\beta(\alpha^2+\dfrac{5\delta^2}{2}+2\alpha\delta+\dfrac{\beta\gamma}{2})\lambda_2\\[2pt]
       \hspace{1cm}+\beta(\alpha^2+\delta^2+\dfrac{3\alpha\delta}{2}+\dfrac{\beta\gamma}{2})\lambda_3=0\\[2pt]
      \beta(-3\alpha^2-2\beta^2+2\delta^2+2\alpha\delta-\beta\gamma)\lambda_1+\delta(-2\alpha^2-2\beta^2+2\delta^2+\alpha\delta)\lambda_3=0\\[2pt]
      \beta(\dfrac{5\alpha^2}{2}+\beta^2+\dfrac{\alpha\delta}{2}+\beta\gamma)\lambda_1+\delta(\alpha^2+\beta^2-\delta^2-\dfrac{\alpha\delta}{2})\lambda_2\\[2pt]
      \hspace{1cm}+\delta(\delta^2+\dfrac{\alpha\delta}{2}+\dfrac{3\beta\gamma}{2})\lambda_3=0\\[2pt]
      -\beta(2\alpha^2+2\delta^2+3\alpha\delta+\beta\gamma)\lambda_1-\delta(2\delta^2+\alpha\delta+3\beta\gamma)\lambda_3=0\\[2pt]
       \end{cases}
\end{eqnarray}
\indent By solving (2.41) , we get
\begin{thm}
$(G_7, g, V)$ admits left-invariant Ricci collineations associated to the Bott connection $\nabla^{B_1}$ if and only if one of the following holds:\\
{\rm (1)}$\alpha= 0,\beta=0,\delta\neq 0, $\\
{\rm (2)}$\alpha= 0,\delta\neq 0,\beta\neq 0,\gamma=0,$\\
{\rm (3)}$\alpha\neq 0,\gamma= 0,\delta= 0.$\\
\indent Moreover, in these cases, we have\\
{\rm (1)}$\mathscr{V}_{\mathscr{R}{C}}=<e_1>$,\\
{\rm (2)}$\mathscr{V}_{\mathscr{R}{C}}=<-\dfrac{\delta}{\beta}e_1+e_2+ e_3>$,\\
{\rm (3)}$\mathscr{V}_{\mathscr{R}{C}}=<e_2+e_3>$.\\
where $\mathscr{V}_{\mathscr{R}{C}}$ is the vector space of left-invariant Ricci collineations on $(G_7, g, V)$.
\end{thm}

\begin{Proof}
We know that $\alpha+\delta\neq 0,\alpha\gamma= 0$.\\
case {\rm 1)} If $\alpha=0$, then $\delta\neq 0$. By the first equation, $\beta^2\delta(\lambda_2-\lambda_3)=0$.\\
case {\rm 1-1)} If $\beta=0$, by {\rm (2.41)}
\begin{eqnarray}
\begin{cases}
\delta^3\lambda_3=0\\[2pt]
-\delta^3\lambda_2+\delta^3\lambda_3=0\\[2pt]
\end{cases}
\end{eqnarray}
then  $\lambda_2=\lambda_3 =0$. We get {\rm (1)}.\\
case {\rm 1-2)} If $\beta\neq 0$, then $\lambda_2=\lambda_3$. By the second equation, we have $3\beta^2\delta\lambda_1+3\beta\delta^2\lambda_2=0$, i.e. $\beta\lambda_1+\delta\lambda_2=0.$ Then by the sixth equation, we have $2\beta\gamma\delta\lambda_2=0$, i.e. $\gamma\lambda_2=0$.\\
case {\rm 1-2-1)} If $\gamma\neq 0$, then $\lambda_2=\lambda_3=0,\lambda_1=0$.\\
case {\rm 1-2-2)} If $\gamma=0$, then $\beta\lambda_1+\delta\lambda_2=0$, i.e. $\lambda_1=-\dfrac{\delta}{\beta}\lambda_2=-\dfrac{\delta}{\beta}\lambda_3$. {\rm (2.41)} trivially holds. We get {\rm (2)}.\\
case {\rm 2)} If $\alpha\neq 0,$ then $\gamma=0,\alpha+\delta\neq 0$.\\
case {\rm 2-1)} If $2\alpha^3-3\beta^2\delta-\alpha\beta^2\neq 0$, then $\lambda_2=\lambda_3$. By {\rm (2.41)}
\begin{eqnarray}
       \begin{cases}
       (-\alpha^3+\dfrac{3\beta^2\delta}{2}-\dfrac{\alpha\beta^2}{2})\lambda_1+\beta(\dfrac{3\delta^2}{2}+\dfrac{\alpha\delta}{2})\lambda_2=0\\[2pt]
      \beta(-3\alpha^2-2\beta^2+2\delta^2+2\alpha\delta)\lambda_1+\delta(-2\alpha^2-2\beta^2+2\delta^2+\alpha\delta)\lambda_2=0\\[2pt]
      \beta(\dfrac{5\alpha^2}{2}+\beta^2+\dfrac{\alpha\delta}{2})\lambda_1+\delta(\alpha^2+\beta^2)\lambda_2=0\\[2pt]
      -\beta(2\alpha^2+2\delta^2+3\alpha\delta)\lambda_1-\delta(2\delta^2+\alpha\delta)\lambda_3=0\\[2pt]
       \end{cases}
\end{eqnarray}
case {\rm 2-1-1)} If $\delta=0$, then $\lambda_1=0$. {\rm (2.43)} trivially holds. We get {\rm (3)}.\\
case {\rm 2-1-2)} If $\delta\neq 0$, then $\lambda_1=\lambda_2=\lambda_3=0$.\\
case {\rm 2-2)} If $2\alpha^3-3\beta^2\delta-\alpha\beta^2= 0$, by {\rm (2.41)}
\begin{eqnarray}
       \begin{cases}
      \beta(\dfrac{3\alpha^2}{2}+\beta^2-\delta^2-\alpha\delta)\lambda_2+\beta(-\dfrac{3\alpha^2}{2}-\beta^2+\dfrac{5\delta^2}{2}+\dfrac{3\alpha\delta}{2})\lambda_3=0\\[2pt]
      -\beta(\alpha^2+\dfrac{5\delta^2}{2}+2\alpha\delta+\dfrac{\beta\gamma}{2})\lambda_2+\beta(\alpha^2+\delta^2+\dfrac{3\alpha\delta}{2}+\dfrac{\beta\gamma}{2})\lambda_3=0\\[2pt]
      \beta(-3\alpha^2-2\beta^2+2\delta^2+2\alpha\delta)\lambda_1+\delta(-2\alpha^2-2\beta^2+2\delta^2+\alpha\delta)\lambda_3=0\\[2pt]
      \beta(\dfrac{5\alpha^2}{2}+\beta^2+\dfrac{\alpha\delta}{2})\lambda_1+\delta(\alpha^2+\beta^2-\delta^2-\dfrac{\alpha\delta}{2})\lambda_2+\delta(\delta^2+\dfrac{\alpha\delta}{2})\lambda_3=0\\[2pt]
      -\beta(2\alpha^2+2\delta^2+3\alpha\delta)\lambda_1-\delta(2\delta^2+\alpha\delta)\lambda_3=0\\[2pt]
       \end{cases}
\end{eqnarray}
If $\beta=0$, then $\alpha=0$, this is a contradiction, then $\beta\neq 0$.\\
case {\rm 2-1-1)} If $\delta=0$, then $2\alpha^2-\beta^2=0,\lambda_1=0,\lambda_2=\lambda_3$. This case is in {\rm (3)}.\\
case {\rm 2-1-2)} If $\delta\neq 0$, then $\lambda_1=\lambda_2=\lambda_3=0$.\\
\end{Proof}

\section{ Invariant Ricci collineations associated to the Bott connection on three-dimensional Lorentzian Unimodular Lie groups with the second distribution}

\vskip 0.5 true cm
\indent Let $M$ be a smooth manifold, and let $TM=span\{e_1,e_2,e_3\}$, then took the frst distribution: $F_2=span\{e_1,e_3\}$ and ${F_2}^\bot=span\{e_2\}$, where ${e_1,e_2,e_3}$ is a pseudo-orthonormal basis, with $e_3$ timelike.
The  Bott connection $\nabla^{B_2}$ is defined as follows: (see \cite{F}, \cite{J}, \cite{RK})
\begin{eqnarray}
\nabla^{B_2}_XY=
       \begin{cases}
        \pi_{F_2}(\nabla^L_XY),~~~&X,Y\in\Gamma^\infty(F_2) \\[2pt]
       \pi_{F_2}([X,Y]),~~~&X\in\Gamma^\infty({F_2}^\bot),Y\in\Gamma^\infty(F_2)\\[2pt]
       \pi_{{F_2}^\bot}([X,Y]),~~~&X\in\Gamma^\infty(F_2),Y\in\Gamma^\infty({F_2}^\bot)\\[2pt]
       \pi_{{F_2}^\bot}(\nabla^L_XY),~~~&X,Y\in\Gamma^\infty({F_2}^\bot)\\[2pt]
       \end{cases}
\end{eqnarray}
where $\pi_{F_2}$ and $\pi_{F_2}^\bot$ are respectively the projection on $F_2$ and ${F_2}^\bot$, $\nabla^L$ is the Levi-Civita connection of $G_i$ .\\

\vskip 0.5 true cm
\noindent{\bf 3.1 Invariant Ricci collineations of $G_1$ associated to the Bott connection $\nabla^{B_2}$}\\
\vskip 0.5 true cm
\begin{lem}(\cite{Wu}) The Ricci tensor of  $(G_1,g)$ associated to the Bott connection $\nabla^{B_2}$ is determined by
\begin{align}
&Ric^{B_2}(e_1,e_1)=\alpha^2-\beta^2,~~~Ric^{B_2}(e_1,e_2)=\frac{\alpha\beta}{2},~~~Ric^{B_2}(e_1,e_3)=-\alpha\beta,\\\notag
&Ric^{B_2}(e_2,e_2)=0,~~~Ric^{B_2}(e_2,e_3)=\frac{\alpha^2}{2},~~~Ric^{B_2}(e_3,e_3)=\beta^2-\alpha^2.
\end{align}
\end{lem}
By (2.5) and Lemma 3.1, we have
\begin{lem}
\begin{align}
&(L_VRic^{B_2})(e_1,e_1)=2\alpha^3\lambda_1+\alpha(\beta^2-2\alpha^2)\lambda_2,\\\notag
&(L_VRic^{B_2})(e_1,e_2)=-\alpha^3\lambda_1-\beta^3\lambda_3,\\\notag
&(L_VRic^{B_2})(e_1,e_3)=\alpha(\alpha^2-\dfrac{\beta^2}{2})\lambda_1-\dfrac{\alpha^2\beta}{2}\lambda_2+\dfrac{\alpha^2\beta}{2}\lambda_3,\\\notag
&(L_VRic^{B_2})(e_2,e_2)=\alpha(\beta^2+\alpha^2)\lambda_3,\\\notag
&(L_VRic^{B_2})(e_2,e_3)=\beta(\dfrac{\alpha^2}{2}+\beta^2)\lambda_1-\dfrac{\alpha}{2}(\alpha^2+\beta^2)\lambda_2-\dfrac{\alpha^3}{2}\lambda_3,\\\notag
&(L_VRic^{B_2})(e_3,e_3)=-\alpha^2\beta\lambda_1+\alpha^3\lambda_2.
\end{align}
\end{lem}
Then, if a left-invariant vector field $V$ is a  Ricci collineation associated to the Bott connection $\nabla^{B_2}$, by Lemma 3.2 and Theorem 2.1, we have the following equations:
\begin{eqnarray}
       \begin{cases}
       2\alpha^3\lambda_1+\alpha(\beta^2-2\alpha^2)\lambda_2=0\\[2pt]
       -\alpha^3\lambda_1-\beta^3\lambda_3=0\\[2pt]
       \alpha(\alpha^2-\dfrac{\beta^2}{2})\lambda_1-\dfrac{\alpha^2\beta}{2}\lambda_2+\dfrac{\alpha^2\beta}{2}\lambda_3=0\\[2pt]
       \alpha(\beta^2+\alpha^2)\lambda_3=0\\[2pt]
       \beta(\dfrac{\alpha^2}{2}+\beta^2)\lambda_1-\dfrac{\alpha}{2}(\alpha^2+\beta^2)\lambda_2-\dfrac{\alpha^3}{2}\lambda_3=0\\[2pt]
       -\alpha^2\beta\lambda_1+\alpha^3\lambda_2=0\\[2pt]
       \end{cases}
\end{eqnarray}
By solving (3.4) , we get
\begin{thm}
$(G_1, g, V)$ does not admit left-invariant Ricci collineations associated to the Bott connection $\nabla^{B_2}$.
\end{thm}
\begin{Proof}
We know that $\alpha\neq 0$. By the fourth equation, and $\alpha^2+\beta^2=0$, we have $\lambda_3=0$. Then by the second equation, we have $\alpha^3\lambda_1=0,$ i.e. $\lambda_1=0$. Finally by the fifth equation, we have $\alpha(\alpha^2+\beta^2)\lambda_2=0$, i.e. $\lambda_2=0$.
So we have no non-trivial solution.
\end{Proof}

\vskip 0.5 true cm
\noindent{\bf 3.2 Invariant Ricci collineations of $G_2$ associated to the Bott connection $\nabla^{B_2}$}\\
\vskip 0.5 true cm
\begin{lem}(\cite{Wu}) The Ricci tensor of  $(G_2,g)$ associated to the Bott connection $\nabla^{B_2}$ is determined by
\begin{align}
&Ric^{B_2}(e_1,e_1)=-(\beta^2+\gamma^2),~~~Ric^{B_2}(e_1,e_2)=0,~~~Ric^{B_2}(e_1,e_3)=0,\\\notag
&Ric^{B_2}(e_2,e_2)=0,~~~Ric^{B_2}(e_2,e_3)=-\frac{\alpha\gamma}{2},~~~Ric^{B_2}(e_3,e_3)=\gamma^2+\alpha\beta.
\end{align}
\end{lem}
By (2.5) and Lemma 3.4, we have
\begin{lem}
\begin{align}
&(L_VRic^{B_2})(e_1,e_1)=0,\\\notag
&(L_VRic^{B_2})(e_1,e_2)=\dfrac{\alpha\beta\gamma}{2}\lambda_2-\alpha(\beta^2+\dfrac{\gamma^2}{2})\lambda_3,\\\notag
&(L_VRic^{B_2})(e_1,e_3)=\gamma^2(\dfrac{\alpha}{2}-\beta)\lambda_2-\gamma(\dfrac{\alpha\beta}{2}+\gamma^2)\lambda_3,\\\notag
&(L_VRic^{B_2})(e_2,e_2)=-\alpha\beta\gamma\lambda_1,\\\notag
&(L_VRic^{B_2})(e_2,e_3)=\beta(\alpha\beta+\gamma^2)\lambda_1,\\\notag
&(L_VRic^{B_2})(e_3,e_3)=\gamma(\alpha\beta+2\gamma^2)\lambda_1.
\end{align}
\end{lem}
\indent Then, if a left-invariant vector field $V$ is a  Ricci collineation associated to the Bott connection $\nabla^{B_2}$, by Lemma 3.5 and Theorem 2.1, we have the following equations:
\begin{eqnarray}
       \begin{cases}
       \dfrac{\alpha\beta\gamma}{2}\lambda_2-\alpha(\beta^2+\dfrac{\gamma^2}{2})\lambda_3=0\\[2pt]
       \gamma^2(\dfrac{\alpha}{2}-\beta)\lambda_2-\gamma(\dfrac{\alpha\beta}{2}+\gamma^2)\lambda_3=0\\[2pt]
       -\alpha\beta\gamma\lambda_1=0\\[2pt]
       \beta(\alpha\beta+\gamma^2)\lambda_1=0\\[2pt]
       \gamma(\alpha\beta+2\gamma^2)\lambda_1=0\\[2pt]
       \end{cases}
\end{eqnarray}
\indent By solving (3.7) , we get
\begin{thm}
$(G_2, g, V)$ admits left-invariant Ricci collineations associated to the Bott connection $\nabla^{B_2}$ if and only if one of the following holds:\\
{\rm (1)} $\alpha= 0,\beta= 0,\gamma\neq 0, $\\
{\rm (2)} $\alpha=0,\beta\neq 0,\gamma\neq 0$,\\
{\rm (3)} $\alpha\neq 0,\beta\neq 0, \alpha-4\beta=0$.\\
\indent Moreover, in these cases, we have\\
{\rm (1)} $\mathscr{V}_{\mathscr{R}{C}}=<e_2>$,\\
{\rm (2)} $\mathscr{V}_{\mathscr{R}{C}}=<-\dfrac{\gamma}{\beta}e_2+e_3>$,\\
{\rm (3)} $\mathscr{V}_{\mathscr{R}{C}}=<(\dfrac{2\beta}{\gamma}+\dfrac{\gamma}{\beta})e_2+e_3>$.\\
where $\mathscr{V}_{\mathscr{R}{C}}$ is the vector space of left-invariant Ricci collineations on $(G_2, g, V)$.
\end{thm}

\begin{Proof}
We know that $\gamma\neq 0$. \\
case {\rm 1)} If $\lambda_1\neq 0$, by the third equation, then $\alpha\beta=0.$ Note that by the fifth equation, we have
\begin{eqnarray}
\gamma(\alpha\beta+2\gamma^2)\lambda_1=2\gamma^3\lambda_1= 0,\notag
\end{eqnarray}
i.e. $\gamma=0$. This is a contradiction.\\
case {\rm 2)} If $\lambda_1= 0$, then the third,fourth,fifith equations trivially holds.\\
case {\rm 2-1)} If $\alpha=0$, then by the second equation, we have $\gamma^2\beta\lambda_2+\gamma^3\lambda_3= 0$, i.e.
\begin{eqnarray}
 \beta\lambda_2+\gamma\lambda_3= 0.\notag
 \end{eqnarray}
case {\rm 2-1-1)} If $\beta=0$, then $\lambda_3= 0$. We get {\rm (1)}.\\
case {\rm 2-1-2)} If $\beta\neq 0$, then $\beta\lambda_2+\gamma\lambda_3= 0$, i.e. $\lambda_2=-\dfrac{\gamma}{\beta}\lambda_3$. We get {\rm (2)}.\\
case {\rm 2-2)} If $\alpha\neq 0,$ by {\rm (3.7)}, we have
\begin{eqnarray}
\begin{cases}
 \dfrac{\beta\gamma}{2}\lambda_2-(\beta^2+\dfrac{\gamma^2}{2})\lambda_3=0\\[2pt]
 \gamma(\dfrac{\alpha}{2}-\beta)\lambda_2-(\dfrac{\alpha\beta}{2}+\gamma^2)\lambda_3=0\\[2pt]
\end{cases}
\end{eqnarray}
case {\rm 2-2-1)} If $\beta=0$,  by {\rm (3.7)}, we have
\begin{eqnarray}
\begin{cases}
-\dfrac{\gamma^2}{2}\lambda_3=0\\[2pt]
\dfrac{\alpha\gamma}{2}\lambda_2-\gamma^2\lambda_3=0\\[2pt]
\end{cases}
\end{eqnarray}
then we have $\lambda_3=0,\lambda_2=0,$ i.e. $\lambda_1=\lambda_2=\lambda_3=0$.\\
case {\rm 2-2-2)} If $\beta\neq 0$, solve {\rm (3.8)} we have
\begin{eqnarray}
(\dfrac{\alpha}{2}-2\beta)(\beta^2+\gamma^2)\lambda_3= 0,\notag
\end{eqnarray}
If $\lambda_3=0,$ then $\lambda_2=0,$ i.e. $\lambda_1=\lambda_2=\lambda_3=0$.\\
If $\alpha-4\beta=0$, then $\lambda_2=(\dfrac{2\beta}{\gamma}+\dfrac{\gamma}{\beta})\lambda_3$. We get {\rm (3)}.\\
\end{Proof}

\vskip 0.5 true cm
\noindent{\bf 3.3 Invariant Ricci collineations of $G_3$ associated to the Bott connection $\nabla^{B_2}$}\\
\vskip 0.5 true cm
\begin{lem}(\cite{Wu}) The Ricci tensor of  $(G_3,g)$ associated to the Bott connection $\nabla^{B_2}$ is determined by
\begin{align}
&Ric^{B_2}(e_1,e_1)=-\beta\gamma,~~~Ric^{B_2}(e_1,e_2)=0,~~~Ric^{B_2}(e_1,e_3)=0,\\\notag
&Ric^{B_2}(e_2,e_2)=0,~~~Ric^{B_2}(e_2,e_3)=0,~~~Ric^{B_2}(e_3,e_3)=\alpha\beta.
\end{align}
\end{lem}
By (2.5) and Lemma 3.7, we have
\begin{lem}
\begin{align}
&(L_VRic^{B_2})(e_1,e_1)=0,~~~(L_VRic^{B_2})(e_1,e_2)=-\alpha\beta\gamma\lambda_3,\\\notag
&(L_VRic^{B_2})(e_1,e_3)=0,~~~(L_VRic^{B_2})(e_2,e_2)=0,\\\notag
&(L_VRic^{B_2})(e_2,e_3)=\alpha\beta\gamma\lambda_1,~~~(L_VRic^{B_2})(e_3,e_3)=0.
\end{align}
\end{lem}
\indent Then, if a left-invariant vector field $V$ is a  Ricci collineation associated to the Bott connection $\nabla^{B_2}$, by Lemma 3.8 and Theorem 2.1, we have the following equations:
\begin{eqnarray}
       \begin{cases}
       -\alpha\beta\gamma\lambda_3=0\\[2pt]
       \alpha\beta\gamma\lambda_1=0\\[2pt]
       \end{cases}
\end{eqnarray}
\indent By solving (3.12) , we get
\begin{thm}
$(G_3, g, V)$ admits left-invariant Ricci collineations associated to the Bott connection $\nabla^{B_2}$ if and only if one of the following holds:\\
{\rm (1)} $\alpha\beta\gamma=0$ ,\\
{\rm (2)} $\alpha\neq 0,\beta\neq 0,\gamma\neq 0$.\\
\indent Moreover, in these cases, we have\\
{\rm (1)} $\mathscr{V}_{\mathscr{R}{C}}=<e_1,e_2,e_3>$,\\
{\rm (2)} $\mathscr{V}_{\mathscr{R}{C}}=<e_2>$.\\
where $\mathscr{V}_{\mathscr{R}{C}}$ is the vector space of left-invariant Ricci collineations on $(G_3, g, V)$.
\end{thm}
\begin{Proof}
case {\rm 1)} If $\alpha\beta\gamma=0$, the above equations trivially holds. We get {\rm (1)}.\\
case {\rm 2)} If $\alpha\beta\gamma\neq 0$, i.e. $\alpha\neq 0,\beta\neq 0,\gamma\neq 0$, then $\lambda_1=\lambda_3=0$. We get {\rm (2)}\\.
\end{Proof}

\vskip 0.5 true cm
\noindent{\bf 3.4 Invariant Ricci collineations of $G_4$ associated to the Bott connection $\nabla^{B_2}$}\\
\vskip 0.5 true cm
\begin{lem}(\cite{Wu}) The Ricci tensor of  $(G_4,g)$ associated to the Bott connection $\nabla^{B_2}$ is determined by
\begin{align}
&Ric^{B_2}(e_1,e_1)=-(\beta-\eta)^2,~~~Ric^{B_2}(e_1,e_2)=0,~~~Ric^{B_2}(e_1,e_3)=0,\\\notag
&Ric^{B_2}(e_2,e_2)=0,~~~Ric^{B_2}(e_2,e_3)=\dfrac{\alpha}{2},~~~Ric^{B_2}(e_1,e_3)=\alpha\beta+1.
\end{align}
\end{lem}
By (2.5) and Lemma 3.10, we have
\begin{lem}
\begin{align}
&(L_VRic^{B_2})(e_1,e_1)=0,\\\notag
&(L_VRic^{B_2})(e_1,e_2)=\dfrac{\alpha}{2}(2\eta-\beta)\lambda_2+\alpha[\dfrac{1}{2}-(\beta-\eta)^2]\lambda_3,\\\notag
&(L_VRic^{B_2})(e_1,e_3)=(-\dfrac{\alpha}{2}+2\eta-\beta+\alpha\eta^2)\lambda_2+(\dfrac{\alpha\beta}{2}+1)\lambda_3,\\\notag
&(L_VRic^{B_2})(e_2,e_2)=-\alpha(2\eta-\beta)\lambda_1,\\\notag
&(L_VRic^{B_2})(e_2,e_3)=-(2\eta-\beta)(\alpha\beta+1)\lambda_1,\\\notag
&(L_VRic^{B_2})(e_3,e_3)=-(\alpha\beta+2)\lambda_1.
\end{align}
\end{lem}
\indent Then, if a left-invariant vector field $V$ is a  Ricci collineation associated to the Bott connection $\nabla^{B_2}$, by Lemma 3.11 and Theorem 2.1, we have the following equations:
\begin{eqnarray}
       \begin{cases}
       \dfrac{\alpha}{2}(2\eta-\beta)\lambda_2+\alpha[\dfrac{1}{2}-(\beta-\eta)^2]\lambda_3=0\\[2pt]
       (-\dfrac{\alpha}{2}+2\eta-\beta+\alpha\eta^2)\lambda_2+(\dfrac{\alpha\beta}{2}+1)\lambda_3=0\\[2pt]
       -\alpha(2\eta-\beta)\lambda_1=0\\[2pt]
       -(2\eta-\beta)(\alpha\beta+1)\lambda_1=0\\[2pt]
       -(\alpha\beta+2)\lambda_1=0\\[2pt]
       \end{cases}
\end{eqnarray}
\indent By solving (3.15) , we get
\begin{thm}
$(G_4, g, V)$ admits left-invariant Ricci collineations associated to the Bott connection $\nabla^{B_2}$ if and only if one of the following holds:\\
{\rm (1)} $\alpha\beta+2= 0,2\eta-\beta=0,\eta=1~{\rm or}-1, $\\
{\rm (2)} $\alpha=0,2\eta-\beta=0,\eta=1~{\rm or}-1,$\\
{\rm (3)} $\alpha= 0,2\eta-\beta\neq 0,\eta=1~{\rm or}-1,$\\
{\rm (4)} $\alpha\neq 0,2\eta-\beta\neq 0,\eta=1~{\rm or}-1,\alpha=-4(2\eta-\beta),$\\
{\rm (5)} $\alpha\neq 0,\beta=\eta,\eta=1~{\rm or}-1.$\\
\indent Moreover, in these cases, we have\\
{\rm (1)} $\mathscr{V}_{\mathscr{R}{C}}=<e_1>$,\\
{\rm (2)} $\mathscr{V}_{\mathscr{R}{C}}=<e_2>$,\\
{\rm (3)} $\mathscr{V}_{\mathscr{R}{C}}=<-\dfrac{1}{2\eta-\beta} e_2+e_3>$,\\
{\rm (4)} $\mathscr{V}_{\mathscr{R}{C}}=<-(\dfrac{\alpha}{2}+4\eta+\dfrac{4}{\alpha})e_2+e_3>$,\\
{\rm (5)} $\mathscr{V}_{\mathscr{R}{C}}=<e_2-\eta e_3>$.\\
where $\mathscr{V}_{\mathscr{R}{C}}$ is the vector space of left-invariant Ricci collineations on $(G_4, g, V)$.
\end{thm}

\begin{Proof}
We know that $\eta=1~{\rm or}-1$.\\
case {\rm 1)} If $\lambda_1\neq 0$, by the sixth equation, $\alpha\beta+2=0$. Then by the fourth equation, $(2\eta-\beta)\lambda_1=0$, i.e. $2\eta-\beta=0$. By {\rm (3.15)}, we have
\begin{eqnarray}
\begin{cases}
-\dfrac{\alpha}{2}\lambda_3=0\\[2pt]
\dfrac{\alpha}{2}\lambda_2-\lambda_3=0\\[2pt]
\end{cases}
\end{eqnarray}
because $\eta=1~({\rm or} -1)$, we have $\alpha=-1~({\rm or} 1)\neq 0$, i.e. $\lambda_2=\lambda_3=0$. We get {\rm (1)}.\\
case {\rm 2)} If $\lambda_1= 0,$ then the third,fourth,fifith equations trivially holds.\\
case {\rm 2-1)} If $\alpha=0$, then by the second equation, we have $(2\eta-\beta)\lambda_2+\lambda_3=0$.\\
case {\rm 2-1-1)} If $2\eta-\beta= 0$, then $\lambda_3=0$. We get {\rm (2)}.\\
case {\rm 2-1-2)} If $2\eta-\beta\neq 0$, then $\lambda_2=-\dfrac{1}{2\eta-\beta}\lambda_3$. We get {\rm (3)}.\\
case {\rm 2-2)} If $\alpha\neq 0,$ we get
\begin{eqnarray}
\begin{cases}
(2\eta-\beta)\lambda_2+[1-2(\beta-\eta)^2]\lambda_3=0\\[2pt]
(\dfrac{\alpha}{2}+2\eta-\beta)\lambda_2+(\dfrac{\alpha\beta}{2}+1)\lambda_3=0\\[2pt]
\end{cases}
\end{eqnarray}
case {\rm 2-2-1)} If $2\eta-\beta=0$, then $\lambda_3=0, \lambda_2=0$, i.e. $\lambda_1=\lambda_2=\lambda_3=0$.\\
case {\rm 2-2-2)} If $2\eta-\beta\neq 0$, we have $[(\dfrac{\alpha}{2}+2\eta-\beta)[1-2(\beta-\eta)^2]-(\dfrac{\alpha\beta}{2}+1)(2\eta-\beta)]\lambda_2=0$, i.e.
\begin{eqnarray}
[\dfrac{\alpha}{2}+2(2\eta-\beta)](\beta-\eta)^2\lambda_2=0.\notag
\end{eqnarray}
 If $\lambda_2=0$, then $\lambda_3=0$, i.e. $\lambda_1=\lambda_2=\lambda_3=0$.\\
If $\alpha=-4(2\eta-\beta)$, then $\lambda_2=-(\dfrac{\alpha}{2}+\dfrac{4}{\alpha}+4\eta)\lambda_3$. We get {\rm (4)}.\\
If $\beta=\eta$, then $\lambda_3=-\eta\lambda_2$. We get {\rm (5)}.

\end{Proof}

\vskip 0.5 true cm
\noindent{\bf 3.5 Invariant Ricci collineations of $G_5$ associated to the Bott connection $\nabla^{B_2}$}\\
\vskip 0.5 true cm
\begin{lem}(\cite{Wu}) The Ricci tensor of  $(G_5,g)$ associated to the Bott connection $\nabla^{B_2}$ is determined by
\begin{align}
&Ric^{B_2}(e_1,e_1)=\alpha^2,~~~Ric^{B_2}(e_1,e_2)=0,~~~Ric^{B_2}(e_1,e_3)=0,\\\notag
&Ric^{B_2}(e_2,e_2)=0,~~~Ric^{B_2}(e_2,e_3)=0,~~~Ric^{B_2}(e_3,e_3)=-(\alpha^2+\beta\gamma).
\end{align}
\end{lem}
By (2.5) and Lemma 3.13, we have
\begin{lem}
\begin{align}
&(L_VRic^{B_2})(e_1,e_1)=-2\alpha^3\lambda_3,~~~(L_VRic^{B_2})(e_1,e_2)=\alpha^2\gamma\lambda_3,\\\notag
&(L_VRic^{B_2})(e_1,e_3)=-\alpha^3\lambda_1-\alpha^2\gamma\lambda_2~~~(L_VRic^{B_2})(e_2,e_2)=0,\\\notag
&(L_VRic^{B_2})(e_2,e_3)=0,~~~(L_VRic^{B_2})(e_3,e_3)=0.
\end{align}
\end{lem}
\indent Then, if a left-invariant vector field $V$ is a  Ricci collineation associated to the Bott connection $\nabla^{B_2}$, by Lemma 3.14 and Theorem 2.1, we have the following equations:
\begin{eqnarray}
       \begin{cases}
       -2\alpha^3\lambda_3=0\\[2pt]
       \alpha^2\gamma\lambda_3=0\\[2pt]
       -\alpha^3\lambda_1-\alpha^2\gamma\lambda_2=0\\[2pt]
       \end{cases}
\end{eqnarray}
\indent By solving (3.20) , we get
\begin{thm}
$(G_5, g, V)$ admits left-invariant Ricci collineations associated to the Bott connection $\nabla^{B_2}$ if and only if one of the following holds:\\
{\rm (1)} $\alpha= 0,\beta=0,\delta\neq 0, $\\
{\rm (2)} $\alpha\neq 0,\gamma=0,\alpha+\delta\neq 0,\beta\delta= 0,$\\
{\rm (3)} $\alpha\neq 0,\gamma\neq 0,\alpha+\delta\neq 0,\alpha\gamma+\beta\delta= 0.$\\
\indent Moreover, in these cases, we have\\
{\rm (1)} $\mathscr{V}_{\mathscr{R}{C}}=<e_1,e_2,e_3>$,\\
{\rm (2)} $\mathscr{V}_{\mathscr{R}{C}}=<e_2>$,\\
{\rm (3)} $\mathscr{V}_{\mathscr{R}{C}}=<-\dfrac{\gamma}{\alpha} e_2+e_3>$.\\
where $\mathscr{V}_{\mathscr{R}{C}}$ is the vector space of left-invariant Ricci collineations on $(G_5, g, V)$.
\end{thm}
\begin{Proof}
We know that $\alpha+\delta\neq 0,\alpha\gamma+\beta\delta= 0$.\\
case {\rm 1)} If $\alpha=0$, then $\delta\neq 0, \beta=0$. {\rm (3.20)} trivially holds. We get {\rm (1)}.\\
case {\rm 2)} If $\alpha\neq 0,$ by the first equation we have $\lambda_3=0$. By {\rm (3.20)}, we have
\begin{eqnarray}
\alpha\lambda_1+\gamma\lambda_2=0.\notag
\end{eqnarray}
case {\rm 2-1)} If $\gamma=0$, then $\lambda_1=0,\beta\delta= 0$. We get {\rm (2)}.\\
case {\rm 2-2)} If $\gamma\neq 0$, then $\lambda_1=-\dfrac{\gamma}{\alpha}\lambda_2$. We get {\rm (3)}.\\
\end{Proof}

\vskip 0.5 true cm
\noindent{\bf 3.6 Invariant Ricci collineations of $G_6$ associated to the Bott connection $\nabla^{B_2}$}\\
\vskip 0.5 true cm
\begin{lem}(\cite{Wu}) The Ricci tensor of  $(G_6,g)$ associated to the Bott connection $\nabla^{B_2}$ is determined by
\begin{align}
&Ric^{B_2}(e_1,e_1)=-(\delta^2+\beta\gamma),~~~Ric^{B_2}(e_1,e_2)=0,~~~Ric^{B_2}(e_1,e_3)=0,\\\notag
&Ric^{B_2}(e_2,e_2)=0,~~~Ric^{B_2}(e_2,e_3)=0,~~~Ric^{B_2}(e_3,e_3)=\delta^2.
\end{align}
\end{lem}
By (2.5) and Lemma 3.16, we have
\begin{lem}
\begin{align}
&(L_VRic^{B_2})(e_1,e_1)=0,~~~(L_VRic^{B_2})(e_1,e_2)=0,\\\notag
&(L_VRic^{B_2})(e_1,e_3)=\beta\delta^2\lambda_2+\delta^3\lambda_3~~~(L_VRic^{B_2})(e_2,e_2)=0,\\\notag
&(L_VRic^{B_2})(e_2,e_3)=-\beta\delta^2\lambda_1,~~~(L_VRic^{B_2})(e_3,e_3)=-2\delta^3\lambda_1.
\end{align}
\end{lem}
\indent Then, if a left-invariant vector field $V$ is a  Ricci collineation associated to the Bott connection $\nabla^{B_2}$, by Lemma 3.17 and Theorem 2.1, we have the following equations:
\begin{eqnarray}
       \begin{cases}
       \beta\delta^2\lambda_2+\delta^3\lambda_3=0\\[2pt]
       -\beta\delta^2\lambda_1=0\\[2pt]
       -2\delta^3\lambda_1=0\\[2pt]
       \end{cases}
\end{eqnarray}
\indent By solving (3.23) , we get
\begin{thm}
$(G_6, g, V)$ admits left-invariant Ricci collineations associated to the Bott connection $\nabla^{B_1}$ if and only if one of the following holds:\\
{\rm (1)}$\delta= 0,\gamma=0,\alpha\neq 0, $\\
{\rm (2)}$\delta\neq 0,\beta=0,\alpha+\delta\neq 0,\alpha\gamma= 0,$\\
{\rm (3)}$\delta\neq 0,\beta\neq 0,\alpha+\delta\neq 0,\alpha\gamma-\beta\delta= 0.$\\
\indent Moreover, in these cases, we have\\
{\rm (1)}$\mathscr{V}_{\mathscr{R}{C}}=<e_1,e_2,e_3>$,\\
{\rm (2)}$\mathscr{V}_{\mathscr{R}{C}}=<e_2>$,\\
{\rm (3)}$\mathscr{V}_{\mathscr{R}{C}}=< e_2-\dfrac{\beta}{\delta} e_3>$.\\
where $\mathscr{V}_{\mathscr{R}{C}}$ is the vector space of left-invariant Ricci collineations on $(G_6, g, V)$.
\end{thm}
\begin{Proof}
We know that $\alpha+\delta\neq 0,\alpha\gamma-\beta\delta= 0$.\\
case {\rm 1)} If $\delta=0$, then $\alpha\neq 0, \gamma=0$. {\rm (3.23)} trivially holds. We get {\rm (1)}.\\
case {\rm 2)} If $\delta\neq 0$, by the third equation we have $\lambda_1=0$. By {\rm (3.23)}, we have
 \begin{eqnarray}
\beta\lambda_2+\delta\lambda_3=0.\notag
\end{eqnarray}
case {\rm 2-1)} If $\beta=0$, then $\lambda_3=0,\alpha\gamma=0$. We get {\rm (2)}.\\
case {\rm 2-2)} If $\beta\neq 0$, then $\lambda_3=-\dfrac{\beta}{\delta}\lambda_2$. We get {\rm (3)}.\\
\end{Proof}

\vskip 0.5 true cm
\noindent{\bf 3.7 Invariant Ricci collineations of $G_7$ associated to the Bott connection $\nabla^{B_2}$}\\
\vskip 0.5 true cm
\begin{lem}(\cite{Wu}) The Ricci tensor of  $(G_7,g)$ associated to the Bott connection $\nabla^{B_2}$ is determined by
\begin{align}
&Ric^{B_2}(e_1,e_1)=\alpha^2,~~~Ric^{B_2}(e_1,e_2)=\beta(\alpha+\delta),~~~Ric^{B_2}(e_1,e_3)=\dfrac{\beta(\delta-\alpha)}{2},\\\notag
&Ric^{B_2}(e_2,e_2)=0,~~~Ric^{B_2}(e_2,e_3)=\delta^2+\dfrac{\beta\gamma+\alpha\delta}{2},~~~Ric^{B_2}(e_3,e_3)=-\alpha^2+\beta^2-\beta\gamma.
\end{align}
\end{lem}
By (2.5) and Lemma 3.19, we have
\begin{lem}
\begin{align}
&(L_VRic^{B_2})(e_1,e_1)=(-2\alpha^3-3\beta^2\delta-\alpha\beta^2)\lambda_2+(2\alpha^3+3\beta^2\delta+\alpha\beta^2)\lambda_3,\\\notag
&(L_VRic^{B_2})(e_1,e_2)=(\alpha^3+\dfrac{3\beta^2\delta}{2}+\dfrac{\alpha\beta^2}{2})\lambda_1-\beta(\alpha^2+\delta^2+\dfrac{3\alpha\delta}{2}+\dfrac{\beta\gamma}{2})\lambda_2\\\notag
&\hspace{3.5cm}+\beta(\alpha^2+\dfrac{5\delta^2}{2}+2\alpha\delta+\dfrac{\beta\gamma}{2})\lambda_3,\\\notag
&(L_VRic^{B_2})(e_1,e_3)=(-\alpha^3+\dfrac{3\beta^2\delta}{2}+\dfrac{\alpha\beta^2}{2})\lambda_1-\beta(-\dfrac{3\alpha^2}{2}+\beta^2+\dfrac{5\delta^2}{2}+\dfrac{3\alpha\delta}{2}-\dfrac{\beta\gamma}{2})\lambda_2\\\notag
&\hspace{3.5cm}+\dfrac{\beta}{2}(-3\alpha^2+2\beta^2+2\delta^2+2\alpha\delta-\beta\gamma)\lambda_3,\\\notag
&(L_VRic^{B_2})(e_2,e_2)=2\beta(\alpha^2+\delta^2+\dfrac{3\alpha\delta}{2}+\dfrac{\beta\gamma}{2})\lambda_1+\delta(2\delta^2+\alpha\delta+3\beta\gamma)\lambda_3,\\\notag
&(L_VRic^{B_2})(e_2,e_3)=\beta(-\dfrac{5\alpha^2}{2}+\beta^2-\dfrac{\alpha\delta}{2}-\beta\gamma)\lambda_1-\delta(\delta^2+\dfrac{\alpha\delta}{2}+\dfrac{3\beta\gamma}{2})\lambda_2,\\\notag
&\hspace{3.5cm}+\delta(-\alpha^2+\beta^2+\delta^2+\dfrac{\alpha\delta}{2})\lambda_3,\\\notag
&(L_VRic^{B_2})(e_3,e_3)=-2\beta(-\dfrac{3\alpha^2}{2}+\beta^2+\delta^2+\alpha\delta-\dfrac{\beta\gamma}{2})\lambda_1-2\delta(-\alpha^2+\beta^2+\delta^2+\dfrac{\alpha\delta}{2})\lambda_3.
\end{align}
\end{lem}
\indent Then, if a left-invariant vector field $V$ is a  Ricci collineation associated to the Bott connection $\nabla^{B_2}$, by Lemma 3.20 and Theorem 2.1, we have the following equations:
\begin{eqnarray}
       \begin{cases}
       (-2\alpha^3-3\beta^2\delta-\alpha\beta^2)\lambda_2+(2\alpha^3+3\beta^2\delta+\alpha\beta^2)\lambda_3=0\\[2pt]
       (\alpha^3+\dfrac{3\beta^2\delta}{2}+\dfrac{\alpha\beta^2}{2})\lambda_1-\beta(\alpha^2+\delta^2+\dfrac{3\alpha\delta}{2}+\dfrac{\beta\gamma}{2})\lambda_2\\[2pt]
       \hspace{1cm}+\beta(\alpha^2+\dfrac{5\delta^2}{2}+2\alpha\delta+\dfrac{\beta\gamma}{2})\lambda_3=0\\[2pt]
       (-\alpha^3+\dfrac{3\beta^2\delta}{2}+\dfrac{\alpha\beta^2}{2})\lambda_1-\beta(-\dfrac{3\alpha^2}{2}+\beta^2+\dfrac{5\delta^2}{2}+\dfrac{3\alpha\delta}{2}-\dfrac{\beta\gamma}{2})\lambda_2\\[2pt]
       \hspace{1cm}+\dfrac{\beta}{2}(-3\alpha^2+2\beta^2+2\delta^2+2\alpha\delta-\beta\gamma)\lambda_3=0\\[2pt]
      2\beta(\alpha^2+\delta^2+\dfrac{3\alpha\delta}{2}+\dfrac{\beta\gamma}{2})\lambda_1+\delta(2\delta^2+\alpha\delta+3\beta\gamma)\lambda_3=0\\[2pt]
      \beta(-\dfrac{5\alpha^2}{2}+\beta^2-\dfrac{\alpha\delta}{2}-\beta\gamma)\lambda_1-\delta(\delta^2+\dfrac{\alpha\delta}{2}+\dfrac{3\beta\gamma}{2})\lambda_2\\[2pt]
      \hspace{1cm}+\delta(-\alpha^2+\beta^2+\delta^2+\dfrac{\alpha\delta}{2})=0\\[2pt]
     -2\beta(-\dfrac{3\alpha^2}{2}+\beta^2+\delta^2+\alpha\delta-\dfrac{\beta\gamma}{2})\lambda_1-2\delta(-\alpha^2+\beta^2+\delta^2+\dfrac{\alpha\delta}{2})\lambda_3=0\\[2pt]
       \end{cases}
\end{eqnarray}
\indent By solving (3.26) , we get
\begin{thm}
$(G_7, g, V)$ admits left-invariant Ricci collineations associated to the Bott connection $\nabla^{B_2}$ if and only if one of the following holds:\\
{\rm (1)}$\alpha= 0,\beta=0,\delta\neq 0, $\\
{\rm (2)}$\alpha\neq 0,\gamma=0,\alpha+\delta\neq 0,\beta\neq 0,$\\
{\rm (3)}$\alpha\neq 0,\gamma\neq 0,\alpha+\delta\neq 0,\alpha\gamma-\beta\delta= 0.$\\
\indent Moreover, in these cases, we have\\
{\rm (1)}$\mathscr{V}_{\mathscr{R}{C}}=<e_1>$,\\
{\rm (2)}$\mathscr{V}_{\mathscr{R}{C}}=<\dfrac{\delta}{\beta} e_1+e_2+e_3>$,\\
{\rm (3)}$\mathscr{V}_{\mathscr{R}{C}}=<e_2+e_3>$.\\
where $\mathscr{V}_{\mathscr{R}{C}}$ is the vector space of left-invariant Ricci collineations on $(G_7, g, V)$.
\end{thm}

\begin{Proof}
We know that $\alpha+\delta\neq 0,\alpha\gamma= 0$.\\
case {\rm 1)} If $\alpha=0$, then $\delta\neq 0$. By the first equation, $\beta^2\delta(\lambda_2-\lambda_3)=0$.\\
case {\rm 1-1)} If $\beta=0$, by {\rm (3.26)}
\begin{eqnarray}
\begin{cases}
\delta^3\lambda_3=0\\[2pt]
-\delta^3\lambda_2+\delta^3\lambda_3=0\\[2pt]
\end{cases}
\end{eqnarray}
then  $\lambda_2=\lambda_3 =0$. We get {\rm (1)}.\\
case {\rm 1-2)} If $\beta\neq 0$, then $\lambda_2=\lambda_3$. By the second equation, we have $3\beta^2\delta\lambda_1-3\beta\delta^2\lambda_2=0$, i.e.
\begin{eqnarray}
\beta\lambda_1-\delta\lambda_2=0.\notag
\end{eqnarray}
Then by the fifth equation, we have $2\gamma\delta\lambda_2=0$, i.e. $\gamma\lambda_2=0$.\\
case {\rm 1-2-1)} If $\gamma\neq 0$, then $\lambda_2=\lambda_3=0,\lambda_1=0$.\\
case {\rm 1-2-2)} If $\gamma=0$, then $\beta\lambda_1-\delta\lambda_2=0. i.e. \lambda_1=\dfrac{\delta}{\beta}\lambda_2=\dfrac{\delta}{\beta}\lambda_3. $ {\rm (3.26)} trivially holds. We get {\rm (2)}.\\
case {\rm 2)} If $\alpha\neq 0,$ then $\gamma=0,\alpha+\delta\neq 0$.\\
case {\rm 2-1)} If $-2\alpha^3-3\beta^2\delta-\alpha\beta^2\neq 0$, then $\lambda_2=\lambda_3$. By {\rm (3.26)}
\begin{eqnarray}
       \begin{cases}
       (\alpha^3+\dfrac{3\beta^2\delta}{2}-\dfrac{\alpha\beta^2}{2})\lambda_1+\beta(\dfrac{3\delta^2}{2}+\dfrac{\alpha\delta}{2})\lambda_2=0\\[2pt]
       \beta(\alpha^2+\delta^2+\dfrac{3\alpha\delta}{2})\lambda_1+\delta(\delta^2+\dfrac{\alpha\delta}{2})\lambda_2=0\\[2pt]
      \beta(-\dfrac{5\alpha^2}{2}+\beta^2-\dfrac{\alpha\delta}{2})\lambda_1+\delta(-\alpha^2+\beta^2)\lambda_2=0\\[2pt]
     \beta(-\dfrac{3\alpha^2}{2}+\beta^2+\delta^2+\alpha\delta)\lambda_1+\delta(-\alpha^2+\beta^2+\delta^2+\dfrac{\alpha\delta}{2})\lambda_2=0\\[2pt]
       \end{cases}
\end{eqnarray}
case {\rm 2-1-1)} If $\delta=0$, then $\lambda_1=0$. {\rm (3.26)} trivially holds. We get {\rm (3)}.\\
case {\rm 2-1-2)} If $\delta\neq 0$, then $\lambda_1=\lambda_2=\lambda_3=0$.\\
case {\rm 2-2)} If $-2\alpha^3-3\beta^2\delta-\alpha\beta^2= 0$, by {\rm (3.25)}
\begin{eqnarray}
       \begin{cases}
      -\beta(\alpha^2+\delta^2+\dfrac{3\alpha\delta}{2})\lambda_2+\beta(\alpha^2+\dfrac{5\delta^2}{2}+2\alpha\delta)\lambda_3=0\\
       -2\alpha^3\lambda_1-\beta(-\dfrac{3\alpha^2}{2}+\beta^2+\dfrac{5\delta^2}{2}+\dfrac{3\alpha\delta}{2})\lambda_2\\
       \hspace{1cm}+\dfrac{\beta}{2}(-3\alpha^2+2\beta^2+2\delta^2+2\alpha\delta)\lambda_3=0\\
      2\beta(\alpha^2+\delta^2+\dfrac{3\alpha\delta}{2})\lambda_1+\delta(2\delta^2+\alpha\delta)\lambda_3=0\\
      \beta(-\dfrac{5\alpha^2}{2}+\beta^2-\dfrac{\alpha\delta}{2})\lambda_1-\delta(\delta^2+\dfrac{\alpha\delta}{2})\lambda_2
      +\delta(-\alpha^2+\beta^2+\delta^2+\dfrac{\alpha\delta}{2})=0\\
      -2\beta(-\dfrac{3\alpha^2}{2}+\beta^2+\delta^2+\alpha\delta)\lambda_1-2\delta(-\alpha^2+\beta^2+\delta^2+\dfrac{\alpha\delta}{2})\lambda_3=0\\
       \end{cases}
\end{eqnarray}
If $\beta=0$, then $\alpha=0$, this is a contradiction, then $\beta\neq 0$.\\
case {\rm 2-2-1)} If $\delta=0$, then $2\alpha^2-\beta^2=0,\lambda_1=0,\lambda_2=\lambda_3$. This case is in {\rm (3)}.\\
case {\rm 2-2-2)} If $\delta\neq 0$, then $\lambda_1=\lambda_2=\lambda_3=0$.\\
\end{Proof}

\section{ Invariant Ricci collineations associated to the Bott connection on three-dimensional Lorentzian Unimodular Lie groups with the third distribution}

\vskip 0.5 true cm
\indent Let $M$ be a smooth manifold, and let $TM=span\{e_1,e_2,e_3\}$, then took the third distribution: $F_3=span\{e_2,e_3\}$ and ${F_3}^\bot=span\{e_1\}$, where ${e_1,e_2,e_3}$ is a pseudo-orthonormal basis, with $e_3$ timelike.
The  Bott connection $\nabla^{B_3}$ is defined as follows: (see \cite{F}, \cite{J}, \cite{RK})
\begin{eqnarray}
\nabla^{B_3}_XY=
       \begin{cases}
        \pi_{F_3}(\nabla^L_XY),~~~&X,Y\in\Gamma^\infty(F_3) \\[2pt]
       \pi_{F_3}([X,Y]),~~~&X\in\Gamma^\infty({F_3}^\bot),Y\in\Gamma^\infty(F_3)\\[2pt]
       \pi_{{F_3}^\bot}([X,Y]),~~~&X\in\Gamma^\infty(F_3),Y\in\Gamma^\infty({F_3}^\bot)\\[2pt]
       \pi_{{F_3}^\bot}(\nabla^L_XY),~~~&X,Y\in\Gamma^\infty({F_3}^\bot)\\[2pt]
       \end{cases}
\end{eqnarray}
where $\pi_{F_3}$ and $\pi_{F_3}^\bot$ are respectively the projection on $F_3$ and ${F_3}^\bot$, $\nabla^L$ is the Levi-Civita connection of $G_i$ .\\

\vskip 0.5 true cm
\noindent{\bf 4.1 Invariant Ricci collineations of $G_1$ associated to the Bott connection $\nabla^{B_3}$}\\
\vskip 0.5 true cm
\begin{lem}(\cite{Wu}) The Ricci tensor of  $(G_1,g)$ associated to the Bott connection $\nabla^{B_3}$ is determined by
\begin{align}
&Ric^{B_3}(e_1,e_1)=0,~~~Ric^{B_3}(e_1,e_2)=\frac{\alpha\beta}{2},~~~Ric^{B_3}(e_1,e_3)=-\frac{\alpha\beta}{2},\\\notag
&Ric^{B_3}(e_2,e_2)=-\beta^2,~~~Ric^{B_3}(e_2,e_3)=0,~~~Ric^{B_3}(e_3,e_3)=\beta^2.
\end{align}
\end{lem}
By (2.5) and Lemma 4.1, we have
\begin{lem}
\begin{align}
&(L_VRic^{B_3})(e_1,e_1)=\alpha\beta^2\lambda_2-\alpha\beta^2\lambda_3,\\\notag
&(L_VRic^{B_3})(e_1,e_2)=-\dfrac{\alpha\beta^2}{2}\lambda_1+\dfrac{\alpha^2\beta}{2}\lambda_2+\beta(\beta^2-\dfrac{\alpha^2}{2})\lambda_3,\\\notag
&(L_VRic^{B_3})(e_1,e_3)=\dfrac{\alpha\beta^2}{2}\lambda_1-\beta(\beta^2+\dfrac{\alpha^2}{2})\lambda_2+\dfrac{\alpha^2\beta}{2}\lambda_3,\\\notag
&(L_VRic^{B_3})(e_2,e_2)=-\alpha^2\beta\lambda_1-\alpha\beta^2\lambda_3,\\\notag
&(L_VRic^{B_3})(e_2,e_3)=\alpha^2\beta\lambda_1+\dfrac{\alpha\beta^2}{2}\lambda_2+\dfrac{\alpha\beta^2}{2}\lambda_3,\\\notag
&(L_VRic^{B_3})(e_3,e_3)=-\alpha^2\beta\lambda_1-\alpha\beta^2\lambda_2.
\end{align}
\end{lem}
Then, if a left-invariant vector field $V$ is a  Ricci collineation associated to the Bott connection $\nabla^{B_1}$, by Lemma 4.2 and Theorem 2.1, we have the following equations:
\begin{eqnarray}
       \begin{cases}
       \alpha\beta^2\lambda_2-\alpha\beta^2\lambda_3=0\\[2pt]
       -\dfrac{\alpha\beta^2}{2}\lambda_1+\dfrac{\alpha^2\beta}{2}\lambda_2+\beta(\beta^2-\dfrac{\alpha^2}{2})\lambda_3=0\\[2pt]
       \dfrac{\alpha\beta^2}{2}\lambda_1-\beta(\beta^2+\dfrac{\alpha^2}{2})\lambda_2+\dfrac{\alpha^2\beta}{2}\lambda_3=0\\[2pt]
       -\alpha^2\beta\lambda_1-\alpha\beta^2\lambda_3=0\\[2pt]
       \alpha^2\beta\lambda_1+\dfrac{\alpha\beta^2}{2}\lambda_2+\dfrac{\alpha\beta^2}{2}\lambda_3=0\\[2pt]
       -\alpha^2\beta\lambda_1-\alpha\beta^2\lambda_2=0\\[2pt]
       \end{cases}
\end{eqnarray}
By solving (4.4) , we get
\begin{thm}
$(G_1, g, V)$ admits left-invariant Ricci collineations associated to the Bott connection $\nabla^{B_3}$ if and only if $\alpha\neq 0,\beta=0$. Moreover, in this case, we have $\mathscr{V}_{\mathscr{R}{C}}=<e_1,e_2,e_3>$, where $\mathscr{V}_{\mathscr{R}{C}}$ is the vector space of left-invariant Ricci collineations on $(G_1, g, V)$.
\end{thm}

\begin{Proof}
We know that $\alpha\neq 0$. By the first equation, we get $\beta^2(\lambda_2-\lambda_3)=0$, so\\
case {\rm 1)} If $\beta=0$, {\rm (4.4)} trivially holds. We get solution.\\
case {\rm 2)} If $\beta\neq 0$, then $\lambda_2=\lambda_3$. By {\rm (4.4)},
\begin{eqnarray}
       \begin{cases}
       -\dfrac{\alpha\beta^2}{2}\lambda_1+\beta^3\lambda_2=0\\[2pt]
       -\alpha^2\beta\lambda_1-\alpha\beta^2\lambda_2=0\\[2pt]
       \end{cases}
\end{eqnarray}
i.e.
\begin{eqnarray}
       \begin{cases}
       \alpha\lambda_1-2\beta\lambda_2=0\\[2pt]
       \alpha\lambda_1+\beta\lambda_2=0\\[2pt]
       \end{cases}
\end{eqnarray}
then we get $\lambda_2=0$, i.e. $\lambda_1=\lambda_2=\lambda_3=0$.\\
\end{Proof}

\vskip 0.5 true cm
\noindent{\bf 4.2 Invariant Ricci collineations of $G_2$ associated to the Bott connection $\nabla^{B_3}$}\\
\vskip 0.5 true cm
\begin{lem}(\cite{Wu}) The Ricci tensor of  $(G_2,g)$ associated to the Bott connection $\nabla^{B_3}$ is determined by
\begin{align}
&Ric^{B_3}(e_1,e_1)=0,~~~Ric^{B_3}(e_1,e_2)=0,~~~Ric^{B_3}(e_1,e_3)=0,\\\notag
&Ric^{B_3}(e_2,e_2)=-\alpha\beta,~~~Ric^{B_3}(e_2,e_3)=-\alpha\gamma,~~~Ric^{B_3}(e_3,e_3)=\alpha\beta.
\end{align}
\end{lem}
By (2.5) and Lemma 4.4, we have
\begin{lem}
\begin{align}
&(L_VRic^{B_3})(e_1,e_1)=0,~~~(L_VRic^{B_3})(e_1,e_2)=\alpha(\beta^2+\gamma^2)\lambda_3,\\\notag
&(L_VRic^{B_3})(e_1,e_3)=-\alpha(\beta^2+\gamma^2)\lambda_2,~~~(L_VRic^{B_3})(e_2,e_2)=0,\\\notag
&(L_VRic^{B_3})(e_2,e_3)=0,~~~(L_VRic^{B_3})(e_3,e_3)=0.
\end{align}
\end{lem}
\indent Then, if a left-invariant vector field $V$ is a  Ricci collineation associated to the Bott connection $\nabla^{B_3}$, by Lemma 4.5 and Theorem 2.1, we have the following equations:
\begin{eqnarray}
       \begin{cases}
       \alpha(\beta^2+\gamma^2)\lambda_3=0\\[2pt]
       -\alpha(\beta^2+\gamma^2)\lambda_2=0\\[2pt]
       \end{cases}
\end{eqnarray}
\indent By solving (4.9) , we get
\begin{thm}
$(G_2, g, V)$ admits left-invariant Ricci collineations associated to the Bott connection $\nabla^{B_3}$ if and only if one of the following holds:\\
{\rm (1)}$\alpha= 0,\gamma\neq 0, $\\
{\rm (2)}$\alpha\neq 0,\gamma\neq 0$.\\
\indent Moreover, in these cases, we have\\
{\rm (1)}$\mathscr{V}_{\mathscr{R}{C}}=<e_1,e_2,e_3>$,\\
{\rm (2)}$\mathscr{V}_{\mathscr{R}{C}}=<e_1>$.\\
where $\mathscr{V}_{\mathscr{R}{C}}$ is the vector space of left-invariant Ricci collineations on $(G_2, g, V)$.
\end{thm}

\begin{Proof}
We know that $\gamma\neq 0$, then $\beta^2+\gamma^2\neq 0$. By {\rm(4.9)}, we have
\begin{eqnarray}
       \begin{cases}
       \alpha\lambda_3=0\\[2pt]
       -\alpha\lambda_2=0\\[2pt]
       \end{cases}
\end{eqnarray}
case {\rm 1)} If $\alpha=0$,  {\rm(4.9)} trivially holds. We get {\rm (1)}.\\
case {\rm 2)} If $\alpha\neq 0$, then $\lambda_2=\lambda_3= 0$. We get {\rm (2)}.\\
\end{Proof}

\vskip 0.5 true cm
\noindent{\bf 4.3 Invariant Ricci collineations of $G_3$ associated to the Bott connection $\nabla^{B_3}$}\\
\vskip 0.5 true cm
\begin{lem}(\cite{Wu}) The Ricci tensor of  $(G_3,g)$ associated to the Bott connection $\nabla^{B_3}$ is determined by
\begin{align}
&Ric^{B_3}(e_1,e_1)=0,~~~Ric^{B_3}(e_1,e_2)=0,~~~Ric^{B_3}(e_1,e_3)=0,\\\notag
&Ric^{B_3}(e_2,e_2)=0,~~~Ric^{B_3}(e_2,e_3)=0,~~~Ric^{B_3}(e_3,e_3)=\alpha\beta.
\end{align}
\end{lem}
By (2.5) and Lemma 4.7, we have
\begin{lem}
\begin{align}
&(L_VRic^{B_3})(e_1,e_1)=0,~~~(L_VRic^{B_3})(e_1,e_2)=0,\\\notag
&(L_VRic^{B_3})(e_1,e_3)=-\alpha\beta\gamma\lambda_2,~~~(L_VRic^{B_3})(e_2,e_2)=0,\\\notag
&(L_VRic^{B_3})(e_2,e_3)=\alpha\beta\gamma\lambda_1,~~~(L_VRic^{B_3})(e_3,e_3)=0.
\end{align}
\end{lem}
\indent Then, if a left-invariant vector field $V$ is a  Ricci collineation associated to the Bott connection $\nabla^{B_3}$, by Lemma 4.8 and Theorem 2.1, we have the following equations:
\begin{eqnarray}
       \begin{cases}
       -\alpha\beta\gamma\lambda_2=0\\[2pt]
       \alpha\beta\gamma\lambda_1=0\\[2pt]
       \end{cases}
\end{eqnarray}
\indent By solving (4.13) , we get
\begin{thm}
$(G_3, g, V)$ admits left-invariant Ricci collineations associated to the Bott connection $\nabla^{B_3}$ if and only if one of the following holds:\\
{\rm (1)}$\alpha\beta\gamma=0$ ,\\
{\rm (2)}$\alpha\neq 0,\beta\neq 0,\gamma\neq 0$.\\
\indent Moreover, in these cases, we have\\
{\rm (1)}$\mathscr{V}_{\mathscr{R}{C}}=<e_1,e_2,e_3>$,\\
{\rm (2)}$\mathscr{V}_{\mathscr{R}{C}}=<e_3>$.\\
where $\mathscr{V}_{\mathscr{R}{C}}$ is the vector space of left-invariant Ricci collineations on $(G_3, g, V)$.
\end{thm}
\begin{Proof}
case {\rm 1)} If $\alpha\beta\gamma=0$, {\rm(4.13)} trivially holds. We get {\rm (1)}.\\
case {\rm 2)} If $\alpha\beta\gamma\neq 0$, i.e. $\alpha\neq 0,\beta\neq 0,\gamma\neq 0$, then $\lambda_1=\lambda_2=0$. We get {\rm (2)}\\.
\end{Proof}

\vskip 0.5 true cm
\noindent{\bf 4.4 Invariant Ricci collineations of $G_4$ associated to the Bott connection $\nabla^{B_3}$}\\
\vskip 0.5 true cm
\begin{lem}(\cite{Wu}) The Ricci tensor of  $(G_4,g)$ associated to the Bott connection $\nabla^{B_3}$ is determined by
\begin{align}
&Ric^{B_3}(e_1,e_1)=0,~~~Ric^{B_3}(e_1,e_2)=0,~~~Ric^{B_3}(e_1,e_3)=0,\\\notag
&Ric^{B_3}(e_2,e_2)=\alpha(2\eta-\beta),~~~Ric^{B_3}(e_2,e_3)=\alpha,~~~Ric^{B_3}(e_3,e_3)=\alpha\beta.
\end{align}
\end{lem}
By (2.5) and Lemma 4.10, we have
\begin{lem}
\begin{align}
&(L_VRic^{B_3})(e_1,e_1)=0,~~~(L_VRic^{B_3})(e_1,e_2)=-\alpha[\beta(2\eta-\beta)-1]\lambda_3,\\\notag
&(L_VRic^{B_3})(e_1,e_3)=\alpha[\beta(2\eta-\beta)-1]\lambda_2,~~~(L_VRic^{B_3})(e_2,e_2)=0,\\\notag
&(L_VRic^{B_3})(e_2,e_3)=0,~~~(L_VRic^{B_3})(e_3,e_3)=0.
\end{align}
\end{lem}
\indent Then, if a left-invariant vector field $V$ is a  Ricci collineation associated to the Bott connection $\nabla^{B_1}$, by Lemma 4.11 and Theorem 2.1, we have the following equations:
\begin{eqnarray}
       \begin{cases}
       -\alpha[\beta(2\eta-\beta)-1]\lambda_3=0\\[2pt]
       \alpha[\beta(2\eta-\beta)-1]\lambda_2=0\\[2pt]
       \end{cases}
\end{eqnarray}
\indent By solving (4.16) , we get
\begin{thm}
$(G_4, g, V)$ admits left-invariant Ricci collineations associated to the Bott connection $\nabla^{B_1}$ if and only if one of the following holds:\\
{\rm (1)}$\alpha=0,\eta=1~{\rm or}-1, $\\
{\rm (2)}$\alpha\neq 0,\beta=\eta,\eta=1~{\rm or}-1,$\\
{\rm (3)}$\alpha\neq 0,\beta\neq\eta,\eta=1~{\rm or}-1.$\\
\indent Moreover, in these cases, we have\\
{\rm (1)}$\mathscr{V}_{\mathscr{R}{C}}=<e_1,e_2,e_3>$,\\
{\rm (2)}$\mathscr{V}_{\mathscr{R}{C}}=<e_1,e_2,e_3>$,\\
{\rm (3)}$\mathscr{V}_{\mathscr{R}{C}}=<e_1>$.\\
where $\mathscr{V}_{\mathscr{R}{C}}$ is the vector space of left-invariant Ricci collineations on $(G_4, g, V)$.
\end{thm}

\begin{Proof}
We know that $\eta=1~{\rm or}-1$.\\
case {\rm 1)} If $\alpha=0$, {\rm(4.16)} trivially holds. We get {\rm (1)}.\\
case {\rm 2)} If $\alpha\neq 0$, then by {\rm(4.16)}, we have
\begin{eqnarray}
       \begin{cases}
      [\beta(2\eta-\beta)-1]\lambda_3=0\\[2pt]
      [\beta(2\eta-\beta)-1]\lambda_2=0\\[2pt]
       \end{cases}
\end{eqnarray}
case {\rm 2-1)} If $\beta(2\eta-\beta)-1=0$, then\\
case {\rm 2-1-1)} If $\eta=1$, we have
\begin{eqnarray}
 \beta(2\eta-\beta)-1=\beta(2-\beta)-1=2\beta-\beta^2-1=-(\beta-1)^2 =0.\notag
\end{eqnarray}
Then $\beta=1$, this case is in {\rm (2)}.\\
case {\rm 2-1-2)} If $\eta=-1$, we have
\begin{eqnarray}
 \beta(2\eta-\beta)-1=\beta(-2-\beta)-1=-2\beta-\beta^2-1=-(\beta+1)^2 =0.\notag
\end{eqnarray}
Then $\beta=-1$, this case is in {\rm (2)}.\\
case {\rm 2-2)} If $\beta(2\eta-\beta)-1\neq 0$, i.e. $\beta\neq\eta$, we have $\lambda_2=\lambda_3=0$. We get {\rm (3)}.\\
\end{Proof}

\vskip 0.5 true cm
\noindent{\bf 4.5 Invariant Ricci collineations of $G_5$ associated to the Bott connection $\nabla^{B_3}$}\\
\vskip 0.5 true cm
\begin{lem}(\cite{Wu}) The Ricci tensor of  $(G_5,g)$ associated to the Bott connection $\nabla^{B_3}$ is determined by
\begin{align}
&Ric^{B_3}(e_1,e_1)=0,~~~Ric^{B_3}(e_1,e_2)=0,~~~Ric^{B_3}(e_1,e_3)=0,\\\notag
&Ric^{B_3}(e_2,e_2)=\delta^2,~~~Ric^{B_3}(e_2,e_3)=0,~~~Ric^{B_3}(e_3,e_3)=-(\delta^2+\beta\gamma).
\end{align}
\end{lem}
By (2.5) and Lemma 4.13, we have
\begin{lem}
\begin{align}
&(L_VRic^{B_3})(e_1,e_1)=0,~~~(L_VRic^{B_3})(e_1,e_2)=\beta\delta^2\lambda_3,\\\notag
&(L_VRic^{B_3})(e_1,e_3)=0,~~~(L_VRic^{B_3})(e_2,e_2)=2\delta^3\lambda_3,\\\notag
&(L_VRic^{B_3})(e_2,e_3)=-\beta\delta^2\lambda_1-\delta^3\lambda_2,~~~(L_VRic^{B_3})(e_3,e_3)=0.
\end{align}
\end{lem}
\indent Then, if a left-invariant vector field $V$ is a  Ricci collineation associated to the Bott connection $\nabla^{B_3}$, by Lemma 4.13 and Theorem 2.1, we have the following equations:
\begin{eqnarray}
       \begin{cases}
       \beta\delta^2\lambda_3=0\\[2pt]
       2\delta^3\lambda_3=0\\[2pt]
       -\beta\delta^2\lambda_1-\delta^3\lambda_2=0\\[2pt]
       \end{cases}
\end{eqnarray}
\indent By solving (4.20) , we get
\begin{thm}
$(G_5, g, V)$ admits left-invariant Ricci collineations associated to the Bott connection $\nabla^{B_3}$ if and only if one of the following holds:\\
{\rm (1)} $\delta= 0,\alpha\neq 0,\gamma=0 $\\
{\rm (2)} $\delta\neq 0,\beta=0,\alpha+\delta\neq 0,\alpha\gamma= 0,$\\
{\rm (3)} $\delta\neq 0,\beta\neq 0,\alpha+\delta\neq 0,\alpha\gamma+\beta\delta= 0.$\\
\indent Moreover, in these cases, we have\\
{\rm (1)} $\mathscr{V}_{\mathscr{R}{C}}=<e_1,e_2,e_3>$,\\
{\rm (2)} $\mathscr{V}_{\mathscr{R}{C}}=<e_1>$,\\
{\rm (3)} $\mathscr{V}_{\mathscr{R}{C}}=<e_1-\dfrac{\beta}{\delta} e_2>$.\\
where $\mathscr{V}_{\mathscr{R}{C}}$ is the vector space of left-invariant Ricci collineations on $(G_5, g, V)$.
\end{thm}
\begin{Proof}
We know that $\alpha+\delta\neq 0,\alpha\gamma+\beta\delta= 0$.\\
case {\rm 1)} If $\delta=0$, then $\alpha\neq 0, \gamma=0$. {\rm (4.20)} trivially holds. We get {\rm (1)}.\\
case {\rm 2)} If $\delta\neq 0$, by the second equation we have $\lambda_3=0$. By {\rm (4.20)}, we have
\begin{eqnarray}
\beta\lambda_1+\delta\lambda_2=0.\notag
\end{eqnarray}
case {\rm 2-1)} If $\beta=0$, then $\lambda_2=0,\alpha\gamma=0$. We get {\rm (2)}.\\
case {\rm 2-2)} If $\beta\neq 0$, then $\lambda_2=-\dfrac{\beta}{\delta}\lambda_1$. We get {\rm (3)}.\\
\end{Proof}

\vskip 0.5 true cm
\noindent{\bf 4.6 Invariant Ricci collineations of $G_6$ associated to the Bott connection $\nabla^{B_3}$}\\
\vskip 0.5 true cm
\begin{lem}(\cite{Wu}) The Ricci tensor of  $(6_6,g)$ associated to the Bott connection $\nabla^{B_3}$ is determined by
\begin{align}
&Ric^{B_1}(e_i,e_j)=0.
\end{align}
for any pairs (i,j).
\end{lem}
By (2.5) and Lemma 4.16, we have
\begin{lem}
\begin{align}
&(L_VRic^{B_1})(e_i,e_j)=0.
\end{align}
\end{lem}
\indent Then we have
\begin{thm}
Any left-invariant vector field on $(G_6, g, V)$ is a left-invariant Ricci collineations associated to the Bott connection $\nabla^{B_3}$ .\\
\end{thm}

\vskip 0.5 true cm
\noindent{\bf 4.7 Invariant Ricci collineations of $G_7$ associated to the Bott connection $\nabla^{B_3}$}\\
\begin{lem}(\cite{Wu}) The Ricci tensor of  $(G_7,g)$ associated to the Bott connection $\nabla^{B_3}$ is determined by
\begin{align}
&Ric^{B_1}(e_1,e_1)=0,~~~Ric^{B_1}(e_1,e_2)=0,~~~Ric^{B_1}(e_1,e_3)=0,\\\notag
&Ric^{B_1}(e_2,e_2)=-\beta\gamma,~~~Ric^{B_1}(e_2,e_3)=\beta\gamma,~~~Ric^{B_3}(e_3,e_3)=-\beta\gamma.
\end{align}
\end{lem}
By (2.5) and Lemma 4.19, we have
\begin{lem}
\begin{align}
&(L_VRic^{B_1})(e_i,e_j)=0.
\end{align}
\end{lem}
\indent Then we have
\begin{thm}
Any left-invariant vector field on $(G_7, g, V)$ is a left-invariant Ricci collineations associated to the Bott connection $\nabla^{B_3}$ .\\
\end{thm}

\section{Acknowledgements}

The author was supported in part by NSFC No.11771070.

\vskip 1 true cm


\bigskip
\bigskip

\noindent {\footnotesize {\it Yanli. Wang} \\
{School of Mathematics and Statistics, Northeast Normal University, Changchun 130024, China}\\
{Email: wangyl553@nenu.edu.cn}

\end{document}